\newif\iffinal
\finaltrue	

\documentclass[letterpaper,11pt,reqno,dvipsnames]{amsart}
\RequirePackage[utf8]{inputenc}
\usepackage[portrait,margin=3cm]{geometry}
\iffinal\else\usepackage[notref,notcite]{showkeys}\fi
\usepackage{mathrsfs,soul}
\usepackage{hyperref}
\usepackage[foot]{amsaddr}
\usepackage{amssymb,amsthm,amsfonts,amsbsy,latexsym,dsfont}
\usepackage[textsize=small]{todonotes}       
\usepackage{graphicx}
\usepackage[numeric,initials,nobysame]{amsrefs}
\usepackage{upref,setspace}
\usepackage{xcolor,colortbl}
\usepackage{enumerate}
\newenvironment{enumeratei}{\begin{enumerate}[\upshape (i)]}{\end{enumerate}}
\newenvironment{enumeratea}{\begin{enumerate}[\upshape (a)]}{\end{enumerate}}

\usepackage{tikz}
\usetikzlibrary{calc,intersections,through,backgrounds,shapes.geometric}
\usetikzlibrary{graphs}
\tikzset{every path/.style={line width=.07 cm}}
 	\usepackage{pgfplots}
\usepackage{paralist}

\newenvironment{inparaenuma}{\begin{inparaenum}[\upshape \bfseries (a) ]}{\end{inparaenum}}

\usepackage{lipsum}

\newcounter{mycounter}

\DeclareMathOperator{\atanh}{atanh}

\numberwithin{equation}{section}
\numberwithin{figure}{section}
\numberwithin{table}{section}

\usepackage{caption}
\usepackage{subcaption}
\usepackage{fourier}
\usepackage{color}

\newtheorem{thm}{Theorem}[section]

\newtheorem{theorem}[thm]{Theorem}

\newtheorem{defn}[thm]{Definition}

\newtheorem{lemma}[thm]{Lemma}

\theoremstyle{definition}

\newtheorem{guess}{Guess}

\renewcommand{\leq}{\leqslant}
\renewcommand{\geq}{\geqslant}

\newcommand{\ind}{\mathds{1}}

\newcommand{\set}[1]{\left\{#1\right\}}

\newcommand{\equald}{\stackrel{\mathrm{d}}{=}}
\newcommand{\probc}{\stackrel{\mathrm{P}}{\longrightarrow}}

\newcommand{\weakc}{\stackrel{\mathrm{w}}{\longrightarrow}}
\newcommand{\convas}{\stackrel{\mathrm{a.e.}}{\longrightarrow}}

\newcommand{\cB}{\mathcal{B}}
\newcommand{\cD}{\mathcal{D}}\newcommand{\cE}{\mathcal{E}}
\newcommand{\cG}{\mathcal{G}}\newcommand{\cH}{\mathcal{H}}
\newcommand{\cL}{\mathcal{L}}
\newcommand{\cM}{\mathcal{M}}\newcommand{\cN}{\mathcal{N}}
\newcommand{\cP}{\mathcal{P}}\newcommand{\cR}{\mathcal{R}}
\newcommand{\cS}{\mathcal{S}}\newcommand{\cT}{\mathcal{T}}\newcommand{\cU}{\mathcal{U}}
\newcommand{\cV}{\mathcal{V}}

\newcommand{\vA}{\mathbf{A}}\newcommand{\vB}{\mathbf{B}}

\newcommand{\vI}{\mathbf{I}}

\newcommand{\vQ}{\mathbf{Q}}

\newcommand{\vW}{\mathbf{W}}
\newcommand{\vY}{\mathbf{Y}}\newcommand{\vZ}{\mathbf{Z}}

\newcommand{\vd}{\mathbf{d}}

\newcommand{\vp}{\mathbf{p}}
\newcommand{\vs}{\mathbf{s}}\newcommand{\vt}{\mathbf{t}}
\newcommand{\vx}{\mathbf{x}}

\newcommand{\mvpi}{\boldsymbol{\pi}}

\newcommand{\fB}{\mathfrak{B}}
\newcommand{\fE}{\mathfrak{E}}
\newcommand{\fI}{\mathfrak{I}}

\newcommand{\fP}{\mathfrak{P}}\newcommand{\fR}{\mathfrak{R}}
\newcommand{\fT}{\mathfrak{T}}

\newcommand{\bC}{\mathbb{C}}

\newcommand{\bG}{\mathbb{G}}
\newcommand{\bL}{\mathbb{L}}

\newcommand{\bQ}{\mathbb{Q}}\newcommand{\bR}{\mathbb{R}}
\newcommand{\bT}{\mathbb{T}}

\newcommand{\bZ}{\mathbb{Z}}

\newcommand{\sF}{\mathscr{F}}

\DeclareMathOperator{\E}{\mathbb{E}}
\DeclareMathOperator{\pr}{\mathbb{P}}

\DeclareMathOperator{\tr}{Tr}

\DeclareMathOperator{\EMP}{EMP}

 \DeclareMathOperator{\dist}{dist}
 
 \DeclareMathOperator{\ERRG}{ERRG}
 
 \DeclareMathOperator{\CM}{CM}
 
 \DeclareMathOperator{\BP}{BP}

\DeclareMathOperator{\pois}{Pois}

\DeclareMathOperator{\model}{model}
\DeclareMathOperator{\supp}{supp}
\DeclareMathOperator{\Bin}{Bin}

\newcommand{\sss}{\scriptscriptstyle}

\newcommand{\erdos}{Erd\H{o}s-R\'enyi }

\newcommand{\bbT}{\mathbb{T}}
\newcommand{\TT}{\mathcal{T}}
\newcommand{\bs}{\mathbf{s}}
\newcommand{\bt}{\mathbf{t}}
\newcommand{\probfr}{\stackrel{\mbox{$\operatorname{a.s.}$-\bf fr}}{\longrightarrow}}
\newcommand{\probcrf}{\stackrel{\mbox{$\operatorname{a.s.}$-\bf efr}}{\longrightarrow}}
\newcommand{\prob}{\mathbb{P}}
\newcommand{\bfomega}{{\boldsymbol \omega}}
\newcommand{\Zbold}{{\mathbb{Z}}}
\newcommand{\LWC}{\stackrel{\mbox{\bf LWC}}{\longrightarrow}}
\newcommand{\probLWC}{\stackrel{\mbox{\bf P-LWC}}{\longrightarrow}}

\usepackage{fourier}

\DeclareMathAlphabet{\mathscrbf}{OMS}{mdugm}{b}{n}

\newcommand{\scrP}{\mathscrbf{P}}

\newcommand{\inpr}[2]{\langle #1, #2 \rangle}

\DeclareMathOperator{\age}{Age}

\usepackage{wrapfig}

\begin{document}

\title[Local weak convergence]{Local weak convergence and its applications}

\date{}
\subjclass[2010]{Primary: 60C05, 05C80. }
\keywords{weak convergence on metric spaces, local weak convergence, continuous time branching processes, random trees, random graphs}

\author[Banerjee]{Sayan Banerjee}
\author[Bhamidi]{Shankar Bhamidi}
\address{Department of Statistics and Operations Research, 304 Hanes Hall, University of North Carolina, Chapel Hill, NC 27599}
\email{sayan,bhamidi,jshen,sethpar@email.unc.edu }
\author[Shen]{Jianan Shen}
\author[Young]{Seth Parker Young}
\maketitle
\begin{abstract}
Motivated in part by understanding average case analysis of fundamental algorithms in computer science, and in part by the wide array of network data available over the last decade, a variety of random graph models, with corresponding processes on these objects, have been proposed over the last few years. The main goal of this paper is to give an overview of local weak convergence, which has emerged as a major technique for understanding large network asymptotics for a wide array of functionals and models. As opposed to a survey, the main goal is to try to explain some of the major concepts and their use to junior researchers in the field and indicate potential resources for further reading. 
\end{abstract}


\section{Introduction}
\label{sec:intro}

Starting graduate courses in probability cover myriad applications including weak convergence of {\bf real-valued} random variables, which one often sees initially via convergence properties of cumulative distribution functions.  However this notion has an equivalent definition via convergence of expectations of functions of the random variables. In many cases, only in a more advanced course does one see see that this final definition easily extends to ``completely general'' metric spaces as beautifully articulated in \cite{parthasarathy2005probability} which then went on to play an important role in a host of applications of modern probability ranging from Billingsley's classic \cite{billingsley2013convergence} to optimal transport \cite{ambrosio2005gradient}. The goal of this specific article is to convey developments of this general theory to the study of asymptotics of large discrete random structures, leading to one specific approach now described as \emph{local} weak convergence. The goal of this paper is not a survey length description of this approach (pointers to beautiful survey length treatments for further reading can be found in Section \ref{sec:conc}). Rather, taking a ``leaf'' from KRP's interview \cite{krp-inter}:

\begin{quote}
	\emph{``{\bf Audience:} So according to you, the best way to learn a new subject is by giving seminars to each other?''\\
	{\bf KRP:} Yes, you can catch hold of somebody who can follow you. Then both of you should discuss and
explain to each other, and the process itself will lead you to something new. Gelfand used to say,
``Read the theorems but not their proofs.''
 You try to construct the proof but that may be too difficult for us. So,
everyday learn at least one new lemma and then try to tell somebody about this lemma. While
discussing, I think something new will always crop up but that is how new ideas are born. Studying by
oneself is boring but if you have four people, and you discuss with each other, go together for tea and
discuss again, then life is more interesting. Ranga Rao, Varadarajan, Varadhan and me, always explained
things to each other.''}
\end{quote}
 
Our goal is to describe the power of ``abstraction'' related to weak convergence on metric spaces, pioneered in \cite{parthasarathy2005probability} in the settings of network models and their applications. The main goal is to convey the joy of mathematical discovery via using some canonical results.   Envisioning an interested student with basic technical knowledge at the level of say a measure theoretic probability course, and interest in discrete random structures, our goal is to explain the motivations driving the mathematical development and (our) intuition behind this fundamental technique. Following Paul Halmos's adage that \emph{``The only way to learn mathematics is to do mathematics
''}, we start with motivating questions, including those that have recently arisen in our own research, in Section \ref{sec:mot}. These questions then motivate the development of the theoretical foundations in Section \ref{sec:defn}. 
Next we convey the ``awe'' behind \emph{some} of the canonical results in this area in Section \ref{sec:power}. Since this paper is in no way a survey,  as a stand alone paper can do little justice to the far reaching applications of this technique, we point the interested reader to further material in the Section \ref{sec:conc}.

\section{Motivating examples and questions}
\label{sec:mot}
We start by describing motivating questions in a number of different areas. 

\subsection{Dynamic social network models }
\label{sec:nodal-att-mot}
Social networks now play an important role in society, for example in the diffusion of information across populations of individuals. This has motivated  the study of network valued data where nodes and/or edges have attributes, which modulate the dynamics of both network evolution, and information flow on the network itself. This has spurred the development of probabilistic network models \cites{Karimi:2018,espin2018towards,espin2022inequality,jordan2013geometric} that incorporate three major ingredients:  (a) heterogeneity in edge creation across groups; (b) dynamic (i.e. time dependent) network evolution and (c) popularity bias. Formulated models are then used to understand domain specific questions, including bias in network sampling, PageRank and degree centrality scores and their impact in network ranking and recommendation algorithms. One basic model is as follows. Fix $\gamma \in [0,1]$. Vertices enter the system sequentially at discrete times $n\geq 1$  starting with a base connected graph $\tilde \cG_0$ at time $n=0$. Write $v_n$ for the vertex that enters at time $n$ and $a(v_n)$ for the corresponding attribute;  every vertex $v_n$ has attribute distribution 
        \begin{equation}
            \label{eqn:726}
            a(v_n) \sim \mvpi, \qquad \text{ independent of } \set{\tilde \cG_s: 0\leq s \leq n-1}. 
        \end{equation}
     The dynamics of construction are recursively defined as: for any $n$ and $v\in \tilde \cG_n$, let $\deg(v,n)$ denote the degree of $v$ at time $n$ \textcolor{black}{(if $\tilde \cG_0 = \{v_0\}$, initialize $\deg(v_0,0)=1$)}. For $n \ge 0$, $v_{n+1}$ attaches to the network via a single outgoing edge. This edge is created via sampling  an existing vertex in $\tilde \cG_n$ to attach to, with probabilities (conditionally on $\tilde \cG_n$ and $a(v_{n+1})$) given by: 
\begin{equation}
\label{eqn:912old}
    \pr(v_{n+1} \leadsto v \,| \,\tilde \cG_n, a(v_{n+1}) = a^{\star}) = \frac{\kappa(a(v), a^\star) [\deg(v,n)]^\gamma}{\sum_{v^\prime \in \tilde \cG_n } \kappa(a(v^\prime), a^\star) [\deg(v^\prime, n)]^\gamma}, \quad v \in \tilde \cG_n.
\end{equation}
\begin{figure}
      \begin{center}
        \includegraphics[scale=.3]{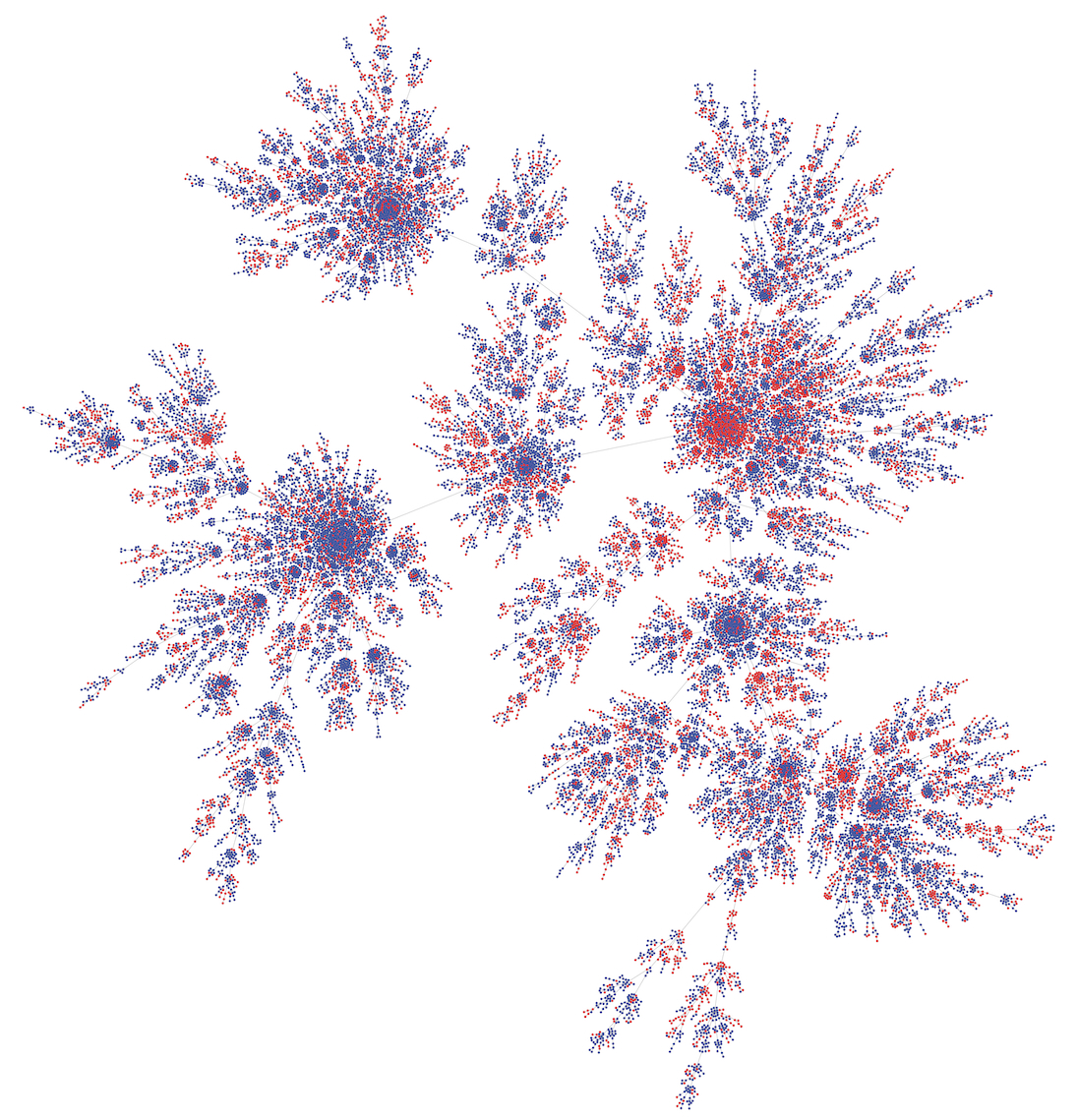}
    \label{fig:nodal}
      \end{center}
    \caption{A simulation with two types ($\cS = \{$red, blue$\}$) with $n=30,000$ and $\kappa(a,a) = .75$, $\kappa(a,a^\prime) = .25$ for $a\neq a^\prime$ and $\pi(red) = .35 $ in the linear $\gamma \equiv 1$ case. }
    \end{figure}

Denote this model of evolving random networks by $\scrP(\gamma, \mvpi, \kappa)$.   The single attribute setting consists of a dynamic network model where new vertices enter the system and connect to existing vertices with probability proportional to their current degree, thus promulgating the ``rich get richer phenomenon'' (i.e. taking degree as a proxy for popularity, high degree vertices accumulate their advantage faster over time). This  was first formulated in the networks community in \cite{barabasi1999emergence}, where they found that this simple mechanistic model gives rise to heavy tailed degree distribution (observed in real world systems); in this specific case, if one writes $N_k(n)$ for the number of vertices with degree $k$ in the network of size $n$, then \cite{barabasi1999emergence} found numerically, which was later proved rigorously in \cite{bollobas2001degree}, that $N_k(n)/n \convas p_k$ as $n\to\infty$, where $p_k\sim C/k^\tau$ as $k\to\infty$ for some $\tau>0$. \\

\subsection{Google's PageRank} PageRank  \cite{page1999pagerank} is one of the most important measures of centrality of vertices in the system. Unlike degree centrality, where the popularity of a vertex is solely quantified by the number of its neighbors, PageRank centrality can also be enhanced by attaching to vertices having high PageRank (popularity enhanced by association with other popular members of a social network). Mathematically, this gives rise to the following definition.

View the networks under consideration as directed trees with edges pointing from offspring to parents. 

    \begin{defn}[PageRank scores with damping factor $c$]
    \label{def:PageRank}
        For a directed graph $\cG = (\cV, \cE)$, the PageRank scores of vertices $v\in \cV$ with damping factor $c$ is the stationary distribution $(\fR_{v,c}: v\in \cG)$ of the following random walk: at each step, with probability $c$, follow an outgoing edge (uniform amongst available choices) from the current location in the graph while, with probability $1-c$, restart at a uniformly selected vertex in the entire graph. These scores are given by the linear system of equations:
        \begin{equation}
            \label{eqn:PageRank}
            \fR_{v,c} = \frac{1-c}{n} + c\sum_{u\in \cN^{-}(v)} \frac{\fR_{u,c}}{d^{+}(u)}
        \end{equation}
        where $\cN^{-}(v)$ is the set of vertices with edges pointed at $v$ and $d^{+}(u)$ is the out-degree of vertex $u$. 
    \end{defn}
    
\noindent {\bf Motivating question:} Understand asymptotics for the empirical distribution of PageRank scores as $n\to\infty$; in particular does this distribution have a limit? What is the relationship, if any, between this limit distribution and the degree distribution? In the attributed network setting, does extremal behavior of the PageRank distribution vary across different attributes? 




\subsection{Gibbs distribution on graphs}
\label{sec:gibbs-distr-intro}
Suppose $\cG= (V(\cG), E(\cG))$ is the vertex and edge representation of a finite connected graph describing a social network. Further suppose, for a specific topic of interest, each individual $v\in \cG$ has an opinion $x_v \in \set{-1,1}$ (e.g. two types of political parties); more complex opinion spaces can be considered but we want to describe an important special case to fix ideas. Now in principle opinions do not necessarily vary independently, rather one would expect opinions of an individual to be influenced at least by it's neighbors in the network. To model such situations, fix two constants $\beta \geq 0, B\in \bR$ and consider the probability measure on opinion configurations $\Omega_{\cG}:= \set{\vx = (x_v:v\in\cG), x_v\in \set{-1,1}~\forall v}$ given by 
\begin{equation}
\label{eqn:gibbs}
\mu(\vx) = \frac{1}{Z(\beta, B)} \exp\left(\beta \sum_{(u,v)\in E(\cG)} x_u x_v + B\sum_{u\in V} x_u\right)
\end{equation}
Here: 
\begin{enumeratea}
 \item $\beta \geq 0$ is sometimes referred to as the inverse temperature and describes the inherent influence of neighbors on one's own opinion; $\beta =0$ implies no influence so everyone's opinions vary independently, while increasing $\beta$ denotes increasing (positive or ferromagnetic, namely opinions of neighbors want to align with each other) influence of the neighbors. 
 \item $B$ is sometimes referred to as the external magnetic field, representing baseline propensities of individuals wanting to adopt either of the two opinions; $B>0$ represents inherent tendency towards $+1$, while $B <0$ represents inherent tendency towards $-1$. 
 \item $Z(\beta, B)$ is the normalizing factor that makes $\mu(\cdot)$ a probability measure (i.e. $\sum_{\vx \in \Omega_{\cG}} \mu(\vx) = 1$). 
 \end{enumeratea} 

\noindent{\bf Motivating question:} Understand properties of the above probability measure. This is hard to do without any other assumptions on the graph, so consider a sequence of graphs $\set{\cG_n:n\geq 1}$ with $\cG_n$ a (random) network on $n$ vertices generated according to some model that incorporates aspects of social networks. For each $n$ consider the Gibbs measure in \eqref{eqn:gibbs} on $\cG_n$ with corresponding normalizing factor $Z_n(\beta, B)$. Understand asymptotics as $n\to\infty$.  

\subsection{Spectral distribution of adjacency matrices}
Random matrix theory now has myriad applications in areas such as statistical physics \cite{mehta2004random}, statistics and machine learning \cites{vershynin2018high,wainwright2019high}. In the context of networks, the adjacency matrix of a graph and its variants play an important role in downstream applications such as classifying nodes of different types (sometimes referred to as community detection) \cites{krzakala2013spectral,abbe2018community,bordenave2015non}. 

 \noindent{\bf Motivating question:} Consider a sequence of random network models $\set{\cG_n:n\geq 1}$ and for $n\geq 1$, let $\vA_n$ denote the adjacency matrix of $\cG_n$ with corresponding eigen-values $\set{\lambda_i^{\sss(n)}:q\leq i\leq n}$. Let $\mu^{\sss(n)}_{\EMP} = n^{-1}\sum_{i=1}^n \delta_{\lambda_i^{\sss(n)}}$ denote the corresponding empirical distribution of eigen-values. Understand asymptotics of $\mu^{\sss(n)}_{\EMP}$ as $n\to\infty$.

\subsection{Probabilistic combinatorial optimization}
\label{sec:prob-opt-mot}
To motivate one stream of local weak convergence (geometric random graphs \cite{aldous-fill-2014}) and also illustrate the breadth of this theory let us now discuss a completely different application. Imagine $n$ faculty who have to be assigned $n$ teaching time slots. Suppose the cost for assigning slot $i$ to faculty $j$ is given by $c_{i,j}$. For a fixed faculty $i$, the best assignment would be $\arg \min_{j\in [n]} c_{i,j}$, however rarely can we guarantee that everyone will get their best choice. Any specific choice of assignment of these slots to faculty can be represented by a permutation $\pi:[n]\to [n]$. Define the total cost of a permutation $C_n(\pi) = \sum_{i\in [n]} c_{i, \pi(i)}. $ and let $\Pi_n$ denote the space of all $n!$ permutations on $[n]$. The goal is to compute $\arg \min_{\pi\in \Pi_n} C_n(\pi)$.  Solving this optimization problem is a canonical example of ``hardness'' from a computational complexity perspective in the sense of worst case upper bounds for algorithms. 


\noindent{\bf Motivating question:} Coupled with computational complexity and algorithmic considerations  that deal with worst case situations, there has been an enormous literature deal with \emph{average case} properties of such questions (see \cite{steele} for a beautiful overview). In the context of the assignment problem consider the following: suppose $\set{c_{i,j}:i,j\in [n]}$ are \emph{i.i.d.} positive random variables (below we will assume exponential rate $1$). Can one understand asymptotics for the cost of the optimal assignment $A_n =\min_{\pi\in \Pi_n} C_n(\pi) $ as $n\to\infty$? After decades of work, in 1987 two statistical physicists Mezard and Parisi (the second of whom won the Nobel prize in Physics in 2022 in part for work related to such questions) via non-rigorous arguments in \cite{mezard1987solution} conjectured that 
\[\E(A_n) \to \zeta(2) = \sum_{i=1}^\infty \frac{1}{i^2} = \frac{\pi^2}{6}, \qquad \text{ as } n\to\infty, \]
where $\zeta(\cdot)$ is the Riemann-Zeta function. Is this not amazing? Our goal is to give an idea of proof of how Aldous in \cite{aldous2001zeta} completed a rigorous proof of this result using local weak convergence techniques. 

\section{Definitions and notation}
\label{sec:defn}

\subsection{Weak convergence: General theory}
In an initial treatment of weak convergence, typically when one specializes to sequences of \emph{real valued} random variables $\set{X_n:n\geq 1} \subseteq \bR$ (say with the goal of describing the Central Limit Theorem), when defining the convergence in distribution of $X_n \weakc X_\infty$ to a limit random variable,  one often phrases weak convergence in terms of the (appropriate) convergence of the cumulative distribution functions (cdfs) $F_{X_n}(x) \to F_{X_\infty}(x)$ for all continuity points $x$ of the limit cdf. 

When dealing with sequences $\set{X_n:n\geq 1} \subseteq \cS$ in a general metric space, the above definition does not extend. However the beautiful theory laid by giants in the field of probability in the early part of the last century (such as Prohorov, Skorohod, the ``Indian school'' including Varadarajan, Ranga Rao, Varadhan and KRP) and described so beautifully in the classics including KRP's \cite{parthasarathy2005probability} and Billingsley's \cite{billingsley2013convergence} describe how one can systematically study weak convergence on general (typically Polish) metric space $\cS$. Let $\cB(\cS)$ denote the associated Borel sigma-field (i.e. generated by open sets in $\cS$) and let $\set{\pr_n:n\geq 1}$ be a sequence of probability measures on the measure space $(\cS, \cB(\cS))$ and let $\pr_\infty$ be another probability measure on the same measure space. 

\begin{defn}
    Let $C_b(\cS)$ denote the space of all continuous and bounded functions on $\cS$. Say that $\pr_n \weakc \pr_\infty$ if for all $f\in C_b(\cS)$, 
    \[\int_{\cS} f(s) d\pr_n(s) \to \int_{\cS} f(s) d\pr_{\infty}(s), \qquad \text{ as } n\to\infty.  \]
\end{defn}

\subsection{Local weak convergence (LWC): Intuition and general framework}
\label{sec:lwc-gfw}
As described in the introduction and motivation, the last decade has witnessed an explosion in the formulation and use of discrete random structures to understand various real world phenomenon. For probabilists, the natural question is to understand asymptotics for such models as the system size $\uparrow \infty$. To fix ideas first consider the following classical network model. 

\begin{defn}[\erdos random graph model]
    Fix $\lambda >0$, $n\geq 1$ and vertex set $[n]:= \set{1,2,\ldots, n}$. Say that a random graph with vertex set $[n]$ has distribution  $\cG_n\sim \ERRG(n,\lambda/n)$ if it is generated as follows: place each of the possible ${n \choose 2}$ edges with probability $\lambda/n$, independently across edges.  
\end{defn}

Now suppose we want to understand asymptotics of this random graph model as $n\to\infty$. Before turning to a rigorous setup let us intuitively try to understand this model. A moment's thought reveals that it is probably hopeless to expect convergence of the entire graph uniformly to a limiting object because there are too many vertices and edges to keep track of.
So let us simplify the question: suppose we sample a vertex $V_n\in [n]$ uniformly at random. What does the geometry of this random graph around this vertex look like? Let us start exploring the graph in a depth first manner around this vertex. The number of friends this vertex has $\cN_n(V_n)\stackrel{d}{=} \Bin(n-1, \lambda/n)$ and by Poisson approximation to the Binomial we get $\cN(V_n) \approx \pois(\lambda)$. The reader should convince themselves that, as we continue the depth first exploration from this vertex, at least ``initially'' in the exploration,  the number of new friends of each explored vertex has approximately $\pois(\lambda)$, more or less independent of the previous steps of the explorations; none of these phrases are rigorous but the above mental picture seems to suggest the following (with the key phrases for the reader in bold):

\begin{guess}
    {\bf Asymptotically $n\to\infty$}, the {\bf local} geometry around a {\bf randomly sampled} vertex should look like a Galton-Watson branching process with offspring distribution $\pois(\lambda)$. Since such a process survives iff $\lambda >1$, this also suggests that there is a {\bf large} connected component iff  $\lambda >1$.
\end{guess}
Thus the natural question: how does one setup the mathematical architecture to make the above guess both rigorous and amenable to answering questions such as the second line? This is the goal of this Section (largely following \cite[Chapter 2]{van2023random}) where we will focus on unweighted graphs. To help the reader, we will discuss special cases of trees and Geometric graphs (i.e. the general edge weighted setting) later.  We will first need some notation. Let $\bG$ denote the space of all locally finite graphs (i.e either finite or infinite graphs such that every vertex has finite, not necessarily uniformly bounded, degree); any such graph can be specified by its vertex set and edge set $\cG = (V(\cG), E(\cG))$. Further, with each edge assumed to have length one,  $\cG$ can be viewed as a metric space on $V(\cG)$ with the usual graph distance $d_{\cG}(\cdot, \cdot)$. For $u\in \cG$ and $r\geq 0$ we let $B_{\cG}(u,r)$ denote the ball of radius $r$ about this vertex. Say that two elements $\cG_1, \cG_2 \in \bG$ are \emph{isomorphic} if there exists (at least one) a bijection $\phi: V(\cG_1) \to V(\cG_2)$ which preserves edges, i.e. $\set{u,v} \ in E(\cG_1)$ iff $\set{\phi(u),\phi(v)} \in E(\cG_2)$. Write this as $\cG_1 \cong \cG_2$. 

A {\bf rooted} graph is a graph with a special vertex $o$, namely a pair $(\cG, o)$ where $\cG \in \bG$ and $o\in V(\cG)$; we will sometimes refer to this as the graph $\cG$ \emph{rooted at} $o$.   Let $\bG_{\star}$ denote the space of all locally finite rooted graphs.  Note that for any $r \geq 0$, we can view the ball of radius $r$ about the root $o$, $B_{\cG}(o,r)$ as an element of $\bG_{\star}$, rooted at $o$. For two elements $\cG_{i,\star} = (\cG_i, o_i) \in \bG_\star$ for $i=1,2$ say that $\cG_{1,\star} \cong \cG_{2,\star} $ iff there exists an isomporphism $\phi$ that preserves the root (i.e. $\phi(o_1) = \phi(o_2)$). 

\begin{defn}[$\bG_{\star}$ as a metric space]
  For two rooted graphs   $\set{\cG_{i,\star} = (\cG_i, o_i): i=1,2}$ define the distance 
  \begin{equation}
      d_{\bG_{\star}}(\cG_{1,\star},\cG_{2,\star}):= \frac{1}{1+R^*}, \qquad R^* =\sup\set{r: B_{\cG_{1,\star}}(o_1, r) \cong B_{\cG_{2,\star}}(o_2, r)},
  \end{equation}
  where the relation $\cong$ denotes root preserving isomorphisms. For the rest of the paper, view $\bG_{\star}$ as a metric space with distance $d_{\bG_{\star}}(\cdot, \cdot)$. 
\end{defn}
Two minor comments are in order: \begin{inparaenuma}
    \item If one was being pedantic, we would now start thinking about $\bG_{\star}$ as equivalence classes of graphs but this leads to far too much overhead. 
    \item It turns out that the above metric makes $\bG_{\star}$ a Polish space. Let $\cB(\bG_\star)$ denote the Borel measure space on $\bG_\star$ generated by this metric.  We are now in a position to talk about weak convergence on this space!
\end{inparaenuma}

\begin{defn}[Local weak convergence and the standard construction]
    \begin{enumeratea}
        \item Let $\set{\pr_n:n\geq 1}$ be a sequence of probability measures on $(\bG_{\star}, \cB(\bG_\star))$ and $\pr_\infty$ be another probability measure on the same measure space. If $~\pr_n \weakc \pr_{\infty}$ as $n\to\infty$ then say that this sequence of measures converges in the local weak convergence sense to $\pr_\infty$ and denote this by $~\pr_n \LWC \pr_\infty$. 
        \item Suppose $\cG$ is a (potentially random element) in $\bG_{<\infty}$ (namely the space of unrooted locally finite graphs) each of finite size. Let $V$ be a vertex selected uniformly at random (and independent of any randomization in the construction of $\cG)$ from $\cG$. Let $\cG[V]$ denote the connected component of $V$ rooted at $V$. This operation of converting an unrooted graph in $\bG_{<\infty}$ to a random element of the space of rooted graphs $\bG_{\star}$ is sometimes referred to as the {\bf standard construction} \cite{aldous-fill-2014}. 
        \item Consider a sequence of (potentially random) graphs $\set{\cG_n: n\geq 1} \subseteq \bG_{<\infty}$ and for each $n$ let $\cG_n(V_n) \in \bG_{\star}$ be the standard construction as in (b) with $V_n$ uniformly selected at random from $\cG_n$ and let $\pr_n$ be the corresponding probability distribution of $\cG_n[V_n]$. Let $\pr_\infty$ be another probability measure on $\bG_{\star}$ and let  $\cG_{\infty, \star} \sim \pr_\infty$. Say that $\cG_n \LWC \cG_{\infty, \star}$ if $~\pr_n \weakc \pr_\infty$ as $n\to\infty$. 
    \end{enumeratea}
\end{defn}
A natural question at this stage is: if this is just weak convergence in a specific setting, what is the rationale for associating a new name to this convergence? Quoting \cite{aldous-fill-2014}:
\begin{quote}
    \emph{``The topology on the metric space $\bG_\star$ turns out to give weak convergence of probability measures on $\bG_\star$ a local character that sharply differs from the traditional weak convergence such as one finds in the weak convergence of scaled random walk to Brownian motion. Weak convergence of measures on $\bG_\star$ never involves any rescaling, and the special role of the neighborhoods $B_{\cG}(o,\cdot)$ means that convergence in $\bG_\star$ only informs us about behavior in the neighborhood of the root.''}
\end{quote}
The power of this methodology and reason for its widespread use is that, despite this local nature of convergence, in a host of examples, a careful analysis can leverage this to the convergence of global functionals. However before proceeding, let us give a more concrete description of what this convergence entails which is the goal of the next Section.

\subsection{Putting the ``local'' in LWC: lifting convergence to the space of probability measures}
\begin{wrapfigure}{l}{0.25\textwidth}
	  \begin{center}
		\includegraphics[width=0.24\textwidth]{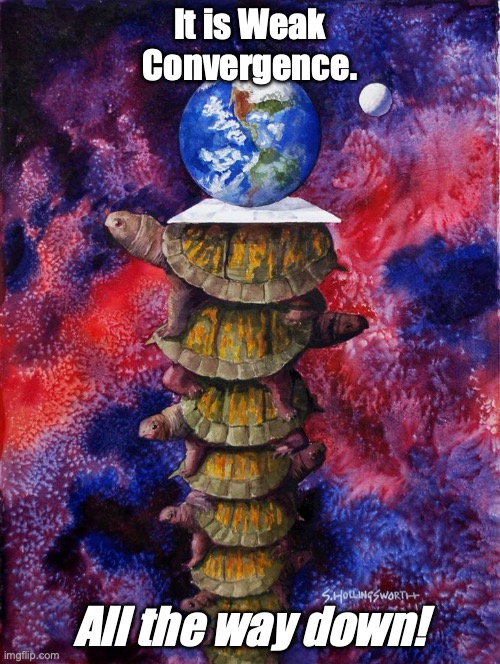}
	  \caption{\tiny Generated via \url{imgflip.com} }
	\label{fig:qst}
	  \end{center}
	\end{wrapfigure} 

This Section perhaps best illustrates the power of weak convergence and the general framework developed in KRP's work on weak convergence and entails unraveling the implications hidden in Definition \ref{def:local-weak}. Let, 
\begin{equation}
    \label{eqn:measure}
    \cM(\bG_{\star}) = \set{\mu: \mu \text{ finite measure on} (\bG_{\star}, \cB(\bG_{\star}))}.
\end{equation}

It turns out, since $\bG_{\star}$ is a Polish space, one can metrize the above space to also make it a Polish space (foundational treatments of such issues can be found in \cite{kallenberg2017random} or \cite[Appendix 2, Vol 1]{daley2003introduction}). Thus now one can talk about convergence of random objects $\set{\fP_n:n\geq 1} \subseteq \cM(\bG_{\star})$!

Now for a finite unrooted graph $\cG \in \bG_{<\infty}$ integer $k\geq 1$ and vertex $v\in \cG$, let $B_{\cG}(v,k) \in \bG_{\star}$ denote the ball of radius $k$ around $v$ in $\cG$, viewed as an element in $\bG_{\star}$ and rooted at $v$. Consider the map:
\begin{equation}
 \fP_k:  \bG_{<\infty} \to \cM(\bG_{\star}) \qquad \text{ given by } \quad  \cG \leadsto \frac{1}{|\cG|} \sum_{v\in \cG} \delta\set{B_{\cG}(v,k)}. 
\end{equation}
If $\pr_\infty$ is a probability measure on $(\bG_{\star}, \cB(\bG_{\star})$, then for $\cG_{\infty,\star} \sim \pr_{\infty}$, write $\fP_{\cG_{\infty,\star},k} = \cL(B_{\cG_{\infty,\star}}(o,k)) $ for the law of the neighborhood upto distance $k$ from the root $o$ of $\cG_{\infty,\star}$.  
\begin{thm}[{\cite[Theorem 2.7]{van2023random}}]
    Let $\set{\cG_n:n\geq 1}$ be a sequence of (potentially random) elements in $\bG_{<\infty}$ and let $\pr_\infty$ be a probability measure on $(\bG_{\star}, \cB(\bG_{\star}))$ and let $\cG_{\infty, \star}\sim \pr_\infty$. Then   $\cG_n \LWC \cG_{\infty, \star}$ iff for every $k\geq 1$ the sequence of (random) empirical neighborhood measures satisfy $\fP_k(\cG_n) \weakc \fP_{\cG_{\infty, \star},k}$ as $n\to\infty$. 
\end{thm}
Thus local weak convergence is equivalent to convergence in distribution of local asymptotics of neighborhoods of randomly sampled vertices. In particular, at first sight, local weak convergence should guarantee convergence of local functionals (e.g. empirical distribution of degrees, empirical distribution of number of vertices at distance two etc) but {\bf not} necessarily global functionals (size of the maximal connected component in the graph etc). 

If we wanted to establish convergence of local statistics in a stronger sense than in distribution (say in probability), the local weak convergence criterion needs to be appropriately strengthened as follows.

\begin{defn}[Local weak convergence in the probability sense]
    \label{def:lwc-prob}
     Let $\set{\cG_n:n\geq 1}$ be a sequence of (potentially random) elements in $\bG_{<\infty}$ and let $\pr_\infty$ be a probability measure on $(\bG_{\star}, \cB(\bG_{\star}))$ and let $\cG_{\infty, \star}\sim \pr_\infty$. Then say that   $\cG_n \probLWC \cG_{\infty, \star}$ written as $\cG_n$ converge in the probability fringe sense to  $\cG_{\infty, \star}$ iff for every $k\geq 1$ the sequence of empirical neighborhood measures satisfy the convergence in probability relationship $\fP_k(\cG_n) \probc \fP_{\cG_{\infty, \star},k}$ as $n\to\infty$. 
\end{defn}

Let us given an example where one has local weak convergence but not convergence in probability. Consider the following sequence of network models $\set{\cG_n:n\geq 1}$ generated as follows:
\begin{enumeratea}
    \item For fixed $n$ flip a biased coin which with probability $1/3$ comes out Heads and $2/3$ Tails. 
    \item If the coin is Head, generate an \erdos random graph $\cG_n \sim \ERRG(n,4/n)$; if it comes out tails let $\cG_n \sim \ERRG(n,10/n)$. 
\end{enumeratea}
Then a reader can convince themselves that, writing $\BP(\pois(\lambda))$ for the genealogy tree generated by Poisson Galton-Watson branching process with mean $\lambda$ started with a single root, that 
\[\cG_n \LWC \frac{1}{3} \BP(\pois(4)) + \frac{2}{3} \BP(\pois(10)),  \]
however this convergence {\bf cannot} be strengthened to convergence in probability. A similar notion is analogously defined for local weak convergence in the almost sure sense. We skip the details.

\subsection{Special case: Local weak convergence of trees}
\label{sec:fringe}

Here we will specialize the phenomenon of local weak convergence to the setting where the sequence of models are rooted trees. There are two main reasons for this: \begin{enumeratea}
    \item This will allow the readers to gain more traction on this notion of convergence in a concrete setting which is easier to visualize (see the figures below) various concepts, than in the context of general graphs. 
    \item More importantly: rooted trees have a notion of direction via orienting the geometry using the root and another specified vertex. When this specified vertex $V_n$ is chosen uniformly at random, it is possible, with some amount of work, to show that for many models of random trees, the subtree {\bf below} this chosen vertex (namely vertices for which the unique path from the root has to pass through this vertex), sometimes called the \emph{fringe distribution}, converges weakly to a distribution on the space of finite trees. What is truly amazing is that, under very general conditions, Aldous in \cite{aldous-fringe} showed that this implies local weak convergence of the tree itself (convergence of entire neighborhoods of the uniformly chosen vertex, not just the subtree below it) to a limiting infinite object! Thus this phenomenon gives a tractable tool to prove local weak convergence as all it requires is to understand what happens ``below'' a randomly selected vertex. This is the main tool we will use in Section \ref{sec:dynamic-res}. 
\end{enumeratea}

 We will first need some notation. For $n\geq 1$, let $ \bT_{n} $ be the space of all rooted trees on  $n$  vertices. Let $ \bbT =
\cup_{n=0}^\infty \bT_{n} $ be the space of all finite rooted marked trees.  Here $\bT_{0} = \emptyset $ will be used to represent the empty tree (tree on zero vertices).  Let $\rho_{\bt}$ denote the root of $\bt$.  For any $r\geq 0$ and $\bt\in \bbT$, let $B(\bt, r) \in \bbT_{\cS}$ denote the subgraph of $\bt$ of vertices within graph distance $r$ from $\rho_{\bt}$, viewed as an element of $\bbT$ and rooted again at $\rho_{\bt}$. 

 Given two rooted finite trees $\bs, \bt \in \bbT$,  say that $\bs \simeq \bt $ if,  there exists a {\bf root
preserving} isomorphism between the two trees viewed as unlabelled graphs.  Given two rooted trees $\bt,\bs \in \bbT$ (see\cite[Equation 2.3.15]{van2023random}), define the distance 
\begin{equation}
\label{eqn:distance-trees}
	d_{\bbT_{\cS}}(\bt,\bs):= \frac{1}{1+R^*}
\end{equation}
where 
\begin{align*}
	R^* =\sup\{r&: B(\bt, r) \simeq B(\bs, r), \text{ and } \exists ~ \text{ isomorphism } { \sF_r} \text{ between }\\
	& B(\bt, r) \text{ and } B(\bs, r) \}.
\end{align*}

\begin{figure}
\centering
\includegraphics[scale=.24]{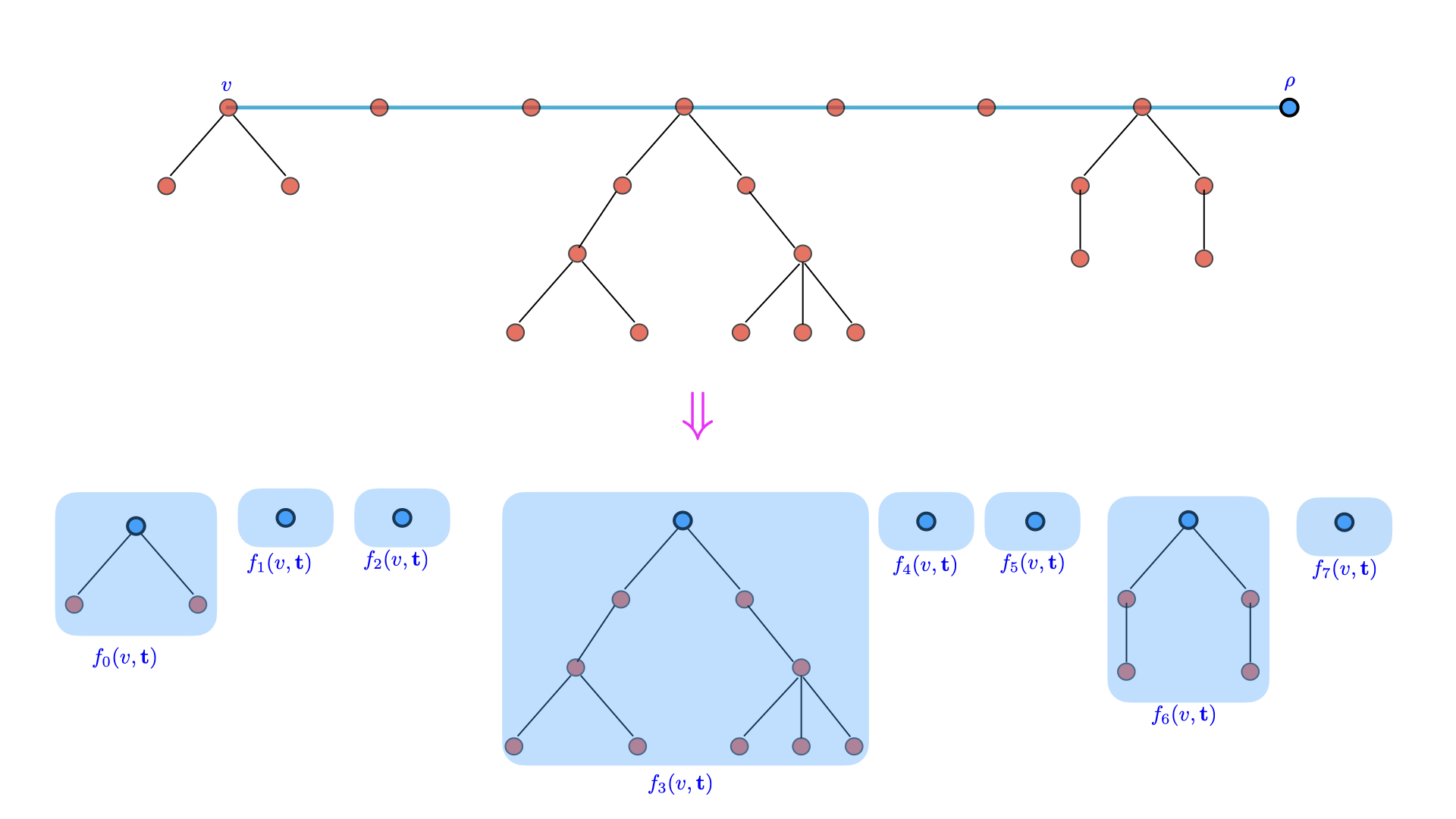}
\caption{Fringe decomposition around vertex $v$ of a finite tree rooted at $\rho$. Here the blue colors represent roots of the respective trees. }
\label{fig:fringe}
\end{figure}

Next, fix a tree $\bt\in \bbT$ with root $\rho = \rho_\bt$ and a vertex $v\in \bt$ at (graph) distance $h$ from the root.  Let $(v_0 =v, v_1, v_2, \ldots, v_h = \rho)$ be the unique path from $v$ to $\rho$. The tree $\bt$  can be decomposed as $h+1$ rooted trees $f_0(v,\bt), \ldots, f_h(v,\bt)$, where $f_0(v,\bt)$ is the tree rooted at $v$, consisting of all vertices for which there exists a path from the root passing through $v$, and for $i \ge 1$, $f_i(v,\bt)$ is the subtree rooted at $v_i$, consisting of all vertices for which the path from the root passes through $v_i$ but not through $v_{i-1}$. Call the map $(v,\bt) \leadsto \bbT^\infty$ where $v\in \bt$, defined via, 
\[F(v, \bt) = \left(f_0(v,\bt), f_1(v,\bt) , \ldots, f_h(v,\bt), \emptyset, \emptyset, \ldots \right),\]
as the fringe decomposition of $\bt$ about the vertex $v$. Call $f_0(v,\bt)$ the {\bf fringe} of the tree $\bt$ at $v$.
For $k\geq 0$, call $F_k(v,\bt) = (f_0(v,\bt) , \ldots, f_k(v,\bt))$ the {\bf extended fringe} of the tree $\bt$ at $v$ truncated at distance $k$ from $v$ on the path to the root ({ see Figure \ref{fig:fringe}}).
Now consider the space $\bbT^\infty$. The metric in \eqref{eqn:distance-trees} can be extended in a straightforward fashion to $\bbT^\infty$.

\begin{figure}
\centering
\includegraphics[scale=.22]{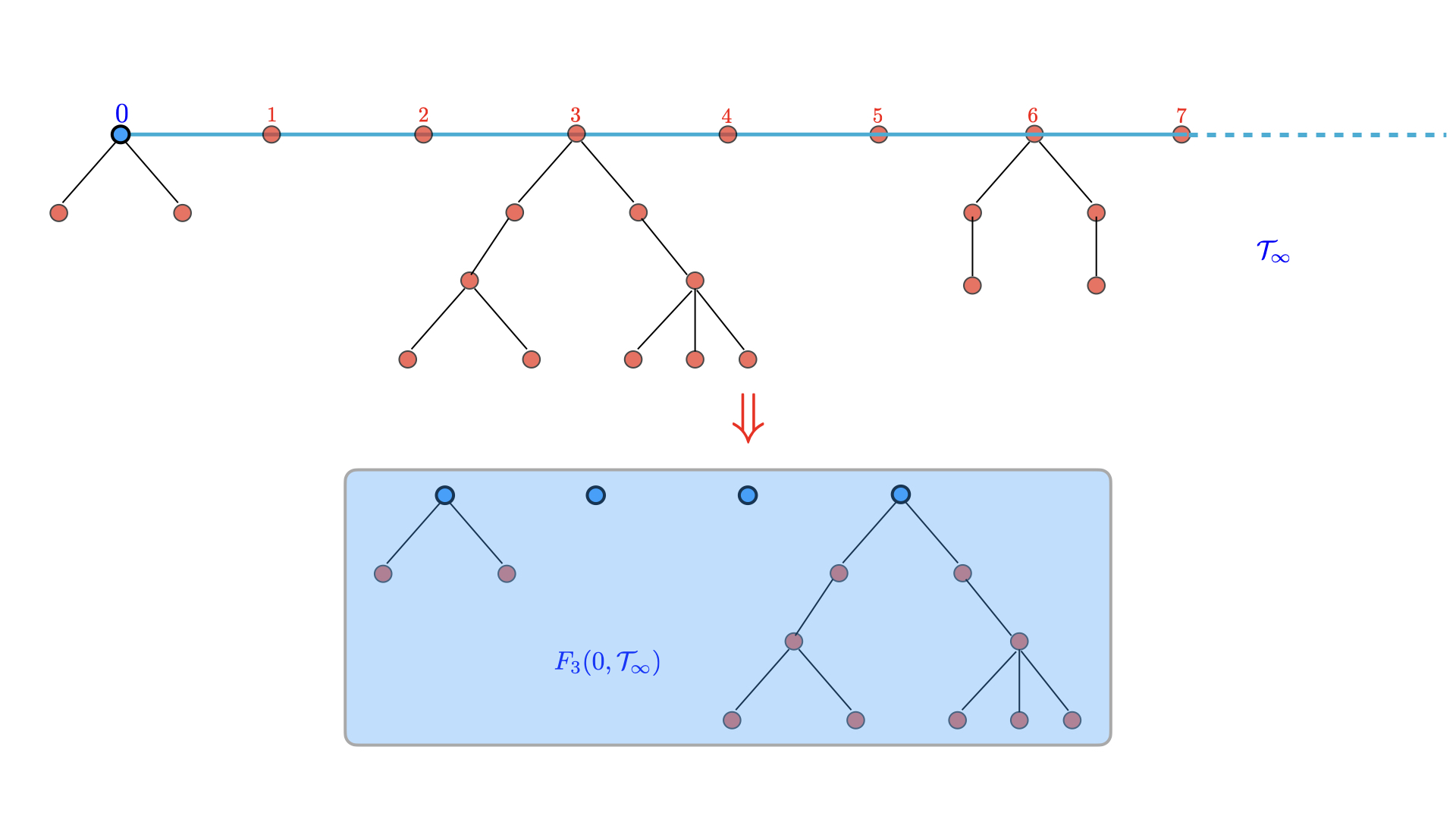}
\caption{A {\tt sin}-tree $\cT_\infty$, namely a tree rooted at $0$ with a single infinite path to infinity, and the corresponding extended fringe $F_3(0,\cT_\infty)$ upto level three about $0$. }
\label{fig:sin}
\end{figure}

Next, an element $\bfomega = (\bt_0, \bt_1, \ldots) \in \bbT^\infty$, with $|\bt_i|\geq 1$ for all $ i\geq 0$,  can be thought of as a locally finite infinite rooted tree with a {\bf s}ingle  path to {\bf in}finity (thus called a {\tt sin}-tree \cite{aldous-fringe}), as follows: identify the sequence of roots of $\set{\bt_i:i\geq 0}$ with the integer lattice $\Zbold_+ = \set{0,1,2,\ldots}$, equipped with the natural nearest neighbor edge set, rooted at $\rho=0$ ({see Figure \ref{fig:sin}}).
  Analogous to the definition of extended fringes for finite trees, for any $k\geq 0$ write 
$F_k(0,\bfomega)= (\bt_0, \bt_1, \ldots, \bt_k)$. 
Call this the extended fringe of the tree $\bfomega$ at vertex $0$, till distance $k$, on the infinite path from $0$. Call $\bt_0 = F_0(0,\bfomega)$ the {\bf fringe} of the {\tt sin}-tree $\bfomega$. Now suppose $\prob$ is a probability measure on $\bbT^\infty$ such that, for $\TT:= (\bt_0(\TT), \bt_1(\TT),\ldots)\sim \prob$,  $|\bt_i(\TT)|\geq 1$ { almost surely (a.s.)}  $\forall~ i\geq 0$. Then $\TT$ can be thought of as an infinite {\bf random} {\tt sin}-tree. 

Define a matrix $\vQ = (\vQ(\vs,\vt): \vs, \vt \in \bbT)$ via the following prescription: suppose the root $\rho_{\vs}$ in $\vs$ has degree $\deg(\rho_{\vs}) \ge 1$, and let $(v_1,\ldots, v_{\deg(\rho_{\vs})})$ denote its children. For $1\leq i\leq v_{\deg(\rho_{\vs})}$, let {$\fT(\vs, v_i)$} be the subtree below $v_i$ and rooted at $v_i$, viewed { as} an element of $\bbT$.
\begin{wrapfigure}{l}{0.25\textwidth}
	  \begin{center}
		\includegraphics[width=0.24\textwidth]{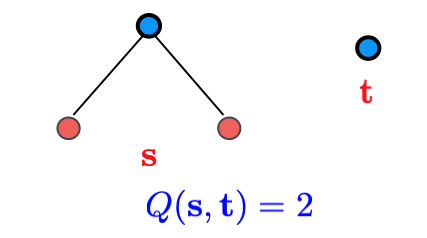}
	\label{fig:qst}
	  \end{center}
	\end{wrapfigure} 
  Write,
\begin{align}
\label{eqn:Q-def}
	\vQ(\vs,\vt):= \sum_{i=1}^{\deg(\rho_{\vs})}\ind\set{d_{\bbT_{\cS}}({ \fT(\vs, v_i)}, \vt) = 0}. 
\end{align} 
See the figure on the left for an example. In words, $Q(\vs, \vt)$ counts the number of descendant subtrees of the root of $\vs$ that are isomorphic to $\vt$. If $\deg(\rho_{\vs})=0$, define $Q(\vs, \vt)=0$. 

Now consider a sequence,
\begin{equation}
\label{eqn:1025}
(\bar \vt_0, \bar \vt_1, \dots) \subseteq \bbT \mbox{ such that } Q(\bar \vt_i, \bar \vt_{i-1}) \ge 1~~ \forall ~ i\geq 1. 
\end{equation}
Then there exists a unique {\tt sin}-tree $\TT$ with infinite path indexed by $\Zbold_+$ such that $\bar \vt_i$ is the subtree rooted at $i$ for all $i \in \Zbold_+$. Conversely, it is easy to see, by taking $\bar \vt_i$ to be the union of (vertices and induced edges) of $\vt_0,\dots, \vt_i$ for each $i \in \Zbold_+$, that every infinite {\tt sin}-tree has a representation of the form \eqref{eqn:1025}. Following \cite{aldous-fringe}, we call this the \emph{monotone representation} of the {\tt sin}-tree $\TT$.

\subsubsection{Convergence on the space of trees}
\label{sec:fringe-convg-def}

Now for any $1\leq k\leq \infty$, let $\cM_{\pr}(\bbT^k)$ denote the space of probability measures on the associated space, metrized using the topology of weak convergence inherited from the corresponding metric on the space $\bbT^k$, see e.g. \cite{billingsley2013convergence}. 
Suppose $\set{\cT_n}_{n\geq 1} \subseteq \bbT$ be a sequence of {\bf finite} rooted random trees on some common probability space (for notational convenience,  assume $|\cT_n|= n$, all one needs is $|\cT_n|\convas \infty$). For $n\geq 1$ and for each fixed $k\geq 0$, consider the empirical distribution of fringes up to distance $k$ ($\delta\{\cdot\}$ below represents Dirac mass):
\begin{equation}
\label{eqn:empirical-fringe-def}
	\fP_{n}^k:= \frac{1}{n} \sum_{v\in \cT_n} \delta\set{F_k(v,\cT_n)}. 
\end{equation}
Thus $\set{\fP_{n}^k:n\geq 1}$ can be viewed as random sequence in  $\cM_{\pr}(\bbT^k)$, with accompanying { notions} of almost sure convergence and convergence in distribution.

\begin{defn}[Local weak convergence]
	\label{def:local-weak}
	Consider two 
	notions of convergence of $\set{\cT_n:n\geq 1}$:
	\begin{enumeratea}
	    \item \label{it:fringe-a} Fix a probability measure $\varpi$ on $\bT$. Say that a sequence of trees $\set{\cT_n}_{n\geq 1}$ converges almost surely, in the fringe sense, to $\varpi$, if 
		\[\fP_n^{0}:= \frac{1}{n} \sum_{v\in \cT_n} \delta\set{F_0(v,\cT_n)} \convas \varpi, \qquad \text{ as } n\to\infty.  \]
	Denote this convergence by $\TT_n\probfr \varpi$ as $n\to\infty$.
	    \item \label{it:fringe-b} Say that a sequence of trees $\set{\cT_n}_{n\geq 1}$ converges almost surely, in the {\bf extended fringe sense}, to a limiting infinite random {\tt sin}-tree $\TT_{\infty}$ if for all $k\geq 0$ one has
	  \[\fP_n^k \stackrel{\mathrm{a.s.}}{\longrightarrow} \prob\left(F_k(0,\TT_{\infty}) \in \cdot \right), \qquad \text{ as } n\to\infty. \]
	Denote this convergence by $\TT_n\probcrf \TT_{\infty}$ as $n\to\infty$.
	\end{enumeratea}
\end{defn} 
{Intuitively, fringe convergence, namely (a) above implies that, if we look at the subtree below \emph{a typical} (i.e. selected uniformly at random) vertex, then the distribution of the corresponding random tree converges in distribution as the network size $n\to\infty$. Extended fringe convergence implies not just the structure of the neighborhood {\bf below} typical vertices, but the entire local neighborhood, within any finite distance converges. Next note that in the setting of the convergence of (b) above, let $\varpi_{\infty}(\cdot) = \pr(F_0(0, \cT_{\infty}) = \cdot)$ denote the distribution of the fringe of $\cT_\infty$ on $\bT_{\cS}$. Convergence in (b) above clearly implies convergence in notion (a) with $\varpi =\varpi_{\infty}(\cdot) $.
More surprisingly, if the limiting distribution $\varpi$ in (a) has a certain `stationarity' property (defined next), \emph{convergence in the fringe sense implies convergence in the extended fringe sense}. }

\begin{defn}[Fringe distribution \cite{aldous-fringe}] \label{fringedef}
	Say that a probability measure $\varpi$ on $\bbT$ is a fringe distribution if 
	\[\sum_{\vs} \varpi(\vs) \vQ(\vs, \vt) = \varpi(\vt), \qquad \forall \ \vt \in \bbT. \]
\end{defn}

For any fringe distribution $\varpi$ on $\bbT_\cS$, one can uniquely obtain the law $\varpi^{EF}$ of a random {\tt sin}-tree $\TT$ with monotone decomposition $(\bar\bt_0(\TT), \bar\bt_1(\TT),\ldots)$ such that for any $i \in \Zbold_+$, any $\bar \vt_0, \bar \vt_1, \dots$ in $\bbT_\cS$,
\begin{equation}\label{ftoef}
\varpi^{EF}((\bar\bt_0(\TT), \bar\bt_1(\TT),\ldots, \bar\bt_i(\TT)) = (\bar \vt_0,\vt_1,\dots,\vt_i)) := \varpi(\vt_i) \prod_{j=1}^{i}Q(\vt_i,\vt_{i-1}),
\end{equation}
where the product is taken to be one if $i=0$. The following Lemma follows by adapting the proof of \cite[Propositions 10 and 11]{aldous-fringe}.
\begin{lemma}\label{ftoeflemma}
Suppose a sequence of trees $\set{\cT_n}_{n\geq 1}$ converges almost surely, in the fringe sense, to $\varpi$. Moreover, suppose that $\varpi$ is a fringe distribution in the sense of Definition \ref{fringedef}. Then $\set{\cT_n}_{n\geq 1}$ converges almost surely, in the extended fringe sense, to a limiting infinite random sin-tree $\TT_{\infty}$ whose law $\varpi^{EF}$ is uniquely obtained from $\varpi$ via \eqref{ftoef}.
\end{lemma}

\textcolor{black}{Both notions imply convergence of functionals such as the degree distribution. For example, in notion (a), letting $\cT_{\varpi} \sim \varpi$ with root denoted by $0$ say, convergence in notion (a) in particular implies,  for any $\ell\geq 0$, 
\begin{equation}
\label{eqn:deg-convg-fr}
	\frac{\#\set{v\in \TT_n, \deg(v) = { \ell} +1}}{n} \convas \prob(\deg(0,\cT_{\varpi})={ \ell}).
\end{equation}
However, both convergences give a lot more information about the asymptotic properties of $\cT_n$, beyond its degree distribution, we will see some specific example in the next Section. } 

\subsection{Extensions}
The basic technical tools outlined above can be extended to various settings. Let us give one example, and leave other settings such as weighted networks (sometimes called Geometric networks) for the reader (and if lost, consult the first few pages of \cite{aldous-steele-obj}). Recall that one of the motivations given in the introduction were attributed network models, namely where vertices have types. Assuming the networks under consideration are trees, let us now show how all the concepts explored in the previous section on fringe convergence easily extend.  Fix attribute space $\cS$ and assume for the rest of the paper that $\cS$ is a Polish space with distance metric $d_{\cS}$. For $n\geq 1$, let $ \bT_{n,\cS} $ be the space of all rooted trees on  $n$  vertices where every vertex has a mark in $\cS$. Let $ \bbT_{\cS} =
\cup_{n=0}^\infty \bT_{n,\cS} $ be the space of all finite rooted marked trees. Here $\bT_{0,\cS} = \emptyset $ will be used to represent the empty tree (tree on zero vertices). For any $\bt\in  \bbT_{\cS} $ and $v\in \bt$, write $a(v)$ for the corresponding attribute of that vertex. Let $\rho_{\bt}$ denote the root of $\bt$.  For any $r\geq 0$ and $\bt\in \bbT_{\cS}$, let $B(\bt, r) \in \bbT_{\cS}$ denote the subgraph of $\bt$ of vertices within graph distance $r$ from $\rho_{\bt}$, viewed as an element of $\bbT_{\cS}$ and rooted again at $\rho_{\bt}$. 

 Given two rooted finite trees $\bs, \bt \in \bbT_{\cS}$,  say that $\bs \simeq \bt $ if, after ignoring all attribute information,  there exists a {\bf root
preserving} isomorphism between the two trees viewed as unlabelled graphs.  Given two rooted trees $\bt,\bs \in \bbT_{\cS}$ (adapting \cite[Equation 2.3.15]{van2023random}), define the distance 
\begin{equation}
\label{eqn:distance-trees2}
	d_{\bbT_{\cS}}(\bt,\bs):= \frac{1}{1+R^*}
\end{equation}
where 
\begin{align*}
	R^* =\sup\{r&: B(\bt, r) \simeq B(\bs, r), \text{ and } \exists ~ \text{ isomorphism } \phi_r \text{ between }\\
	& B(\bt, r) \text{ and } B(\bs, r) \text{ with } d_{\cS}(a(v), a(\phi_r(v))) \leq 1/r~ \forall~ v\in B(\bt, r)
	\}.
\end{align*}
Now proceed with all the notions of convergence, local weak convergence, fringe convergence etc as in the previous Subsection.  

\subsection{Examples}

While all of the results described are well known in the literature, we will give references where a passionate reader new to the field might find complete proofs. 

\subsubsection{Growing random trees}
\label{sec:growing-trees-fringe}
Let us start with an example and describe how the lack of memory property of the exponential distribution transforms this into a more tractable model and then describe the general setting. Consider a sequence of growing random trees $\set{\cT_k:k\geq 1}$, initialized at $k=1$ with $\cT_1$ having one vertex, which we will call the root. For $k\geq 1$, having constructed $\cT_{k-1}$, $\cT_k$ is constructed as follows:  a
  new vertex labelled $k$ enters the system and connects to a randomly selected vertex uniformly at random amongst the available $k-1$ vertices in $\cT_{k-1}$.  Suppose we wanted to understand the asymptotics of this model. 

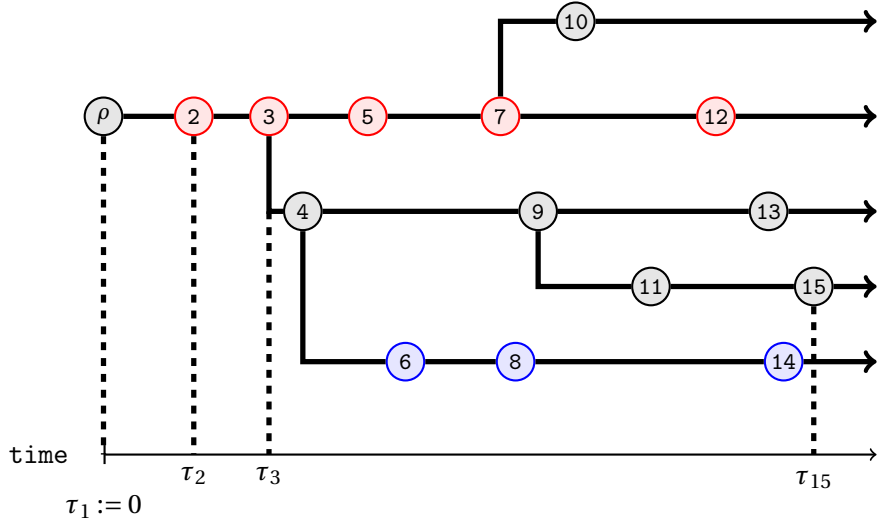
\begin{figure}[h]
\centering
\begin{tikzpicture}[nodes = {draw, circle, thick, align=center, font=\ttfamily\footnotesize, inner sep = 1pt, minimum size = 5mm}, on grid]

 \node [circle, fill=black!10] (rho) {$\rho$};
 \node [circle, fill=red!10, draw=red, right=12mm of rho] (2) {2};
 \node [circle, fill=red!10, draw=red, right=10mm of 2] (3) {3};
 \node [coordinate, right=4.5mm of 3] (4g) { };
 \node [circle, fill=red!10, draw=red, right=6mm of 4g] (5) {5};
 \node [coordinate, right=5mm of 5] (6g) { };
 \node [circle, fill=red!10, draw=red, right=10mm of 6g] (7) {7};
 \node [coordinate, right=2mm of 7] (8g) { };
 \node [coordinate, right=3mm of 8g] (9g) { };
 \node [coordinate, right=5mm of 9g] (10g) { };
 \node [coordinate, right=10mm of 10g] (11g) { };
 \node [circle, fill=red!10, draw=red, right=6mm of 11g] (12) {12};
 \node [coordinate, right=7mm of 12] (13g) { };
 \node [coordinate, right=2mm of 13g] (14g) { };
 \node [coordinate, right=4mm of 14g] (15g) { };

 \node [circle, fill=black!10, below=1cm of 4g] (4) {4};
 \node [circle, fill=blue!10, draw=blue, below=3cm of 6g] (6) {6};
 \node [circle, fill=blue!10, draw=blue, below=3cm of 8g] (8) {8};
 \node [circle, fill=black!10, below=1cm of 9g] (9) {9};
 \node [circle, fill=black!10, above=1cm of 10g] (10) {10};
 \node [circle, fill=black!10, below=2cm of 11g] (11) {11};
 \node [circle, fill=black!10, below=1cm of 13g] (13) {13};
 \node [circle, fill=blue!10, draw=blue, below=3cm of 14g] (14) {14};
 \node [circle, fill=black!10, below=2cm of 15g] (15) {15};

 \path[draw]
 		(rho) -- (2)
 		(2) -- (3)
 		(3) -- (5)
 		(5) -- (7)
 		(7) -- (12)
 		(3) |- (4)
 		(4) -- (9)
 		(9) -- (13)
 		(9) |- (11)
 		(11) -- (15)
 		(4) |- (6)
 		(6) -- (8)
		(8) -- (14)
		(7) |- (10)
 		;

		\draw[->] (10) -- ++ (4,0) coordinate (edge);
		\draw[->] (12) -- (12-|edge);
		\draw[->] (13) -- (13-|edge);
		\draw[->] (14) -- (14-|edge);
		\draw[->] (15) -- (15-|edge);

		\node [coordinate, below=4.5cm of rho, label=left:{\normalsize time$\quad$}] (time0) { };
		\draw[|->, thick] (time0) -- (time0-|edge);

		\draw[dashed] (rho) -- (rho|-time0);
		\draw[dashed] (2) -- (2|-time0);
		\draw[dashed] (3|-4) -- (3|-time0);
		\draw[dashed] (15) -- (15|-time0);

		\node [black, below=1mm, draw=none] at (rho|-time0) {\normalsize$\tau_1 := 0$};
		\node [black, below, draw=none] at (2|-time0) {\normalsize$\tau_2$};
		\node [black, below, draw=none] at (3|-time0) {\normalsize$\tau_3$};
		\node [black, below, draw=none] at (15|-time0) {\normalsize$\tau_{15}$};
\end{tikzpicture}
\caption{The process $\BP_{\alpha}(\cdot)$ in continuous time starting from the root $\rho$ and stopped at $\tau_{15}$.}
\label{fig:test1}
\end{figure}

\begin{figure}[h]
\centering
\begin{tikzpicture}[nodes = {draw, circle, thick, align=center, font=\ttfamily\footnotesize, inner sep = 1pt, minimum size = 5mm}, level 1/.style={sibling distance=4em}, level 2/.style={sibling distance=6em}, level 3/.style={sibling distance=2.5em}]
 \node[fill=black!10] {$\rho$}
   child { node[fill=red!10, draw=red] {2}
     child { node[fill=black!10] {4}
       child { node[fill=blue!10, draw=blue] {6} }
       child { node[fill=blue!10, draw=blue] {8} }
       child { node[fill=blue!10, draw=blue] {14} } }
     child { node[fill=black!10] {9}
       child { node[fill=black!10] {11} }
       child { node[fill=black!10] {15} } }
     child { node[fill=black!10] {13} } }
   child { node[fill=red!10, draw=red] {3} }
   child { node[fill=red!10, draw=red] {5} }
   child { node[fill=red!10, draw=red] {7}
     child { node[fill=black!10] {10} } }
   child { node[fill=red!10, draw=red] {12} };
\end{tikzpicture}
\caption{The corresponding discrete tree containing only the genealogical information of vertices in $\BP_{\alpha}(\tau_{15})$.}
\label{fig:test2}
\end{figure}
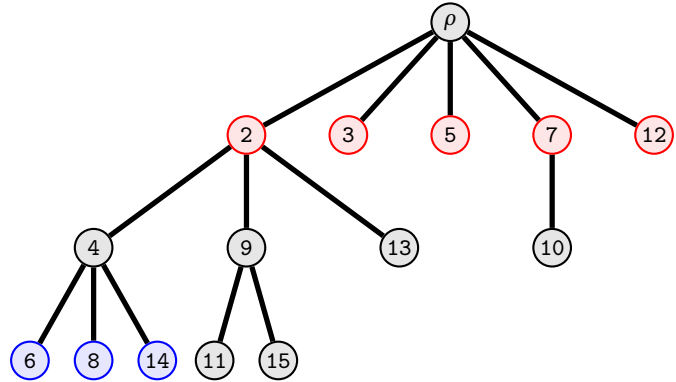

  \noindent{\bf The Zen of the exponential distribution:} We will now consider an embedding of the above tree process in continuous time, an example of the famous Athreya-Karlin embedding \cite{athreya1968embedding} which uses the \emph{lack of memory property} of the exponential distribution. Consider the following \emph{continuous time} branching process $\set{\BP(t):t\geq 0}$ started with a single individual $\rho$ at time $t=0$ with dynamics:
  \begin{enumeratea}
      \item Each individual in the system lives forever. 
      \item The offspring distribution of very individual is a Poisson rate one process independent across individuals; rephrasing, each individual in the branching process gives birth to new individuals at rate one.
  \end{enumeratea}
It turns out, this is a classical and fundamental example of a pure birth process often called the \emph{rate one} Yule process.  For any fixed $t\geq $, let $|\BP(t)|$ denote the size of the branching process at this time and for any $k\geq 1$ define the stopping time $T_k = \inf\set{t: |\BP(t)| = k}$. The following fundamental facts are known about the Yule process and its connection to the discrete time tree process described before:
\begin{enumeratea}
    \item Using the lack of memory property of the exponential distribution, for any $k\geq 1$, consider the  tree $\tilde{\cT}_k$ describing parent child relationships in $\BP(T_k)$ (see Figures \ref{fig:test1} and \ref{fig:test2}.  Then $\tilde{\cT}_k \stackrel{d}{=} {\cT}_k$ and in fact, viewed as growing tree processes, $\set{\tilde{\cT}_k: k\geq 1} \stackrel{d}{=} \set{{\cT}_k: k\geq 1} $. This is one example of the Athreya-Karlin embedding of discrete time processes in continuous time processes. The key conceptual point is that, in the dynamics of the continuous time process, there is a \emph{plethora of independence} since every vertex behaves independently and so is much more amenable to probabilistic analysis. 
    \item\label{it:exp} The process $|\BP(t)|$ grows exponentially at rate one in the sense that there exists a rate one exponential random variable $W$ such that $e^{-t}|\BP(t)| \to W$ a.s. and in $\bL^2$.
\end{enumeratea}
\noindent{\bf Rationale for fringe convergence of $\set{{\cT}_k: k\geq 1}$:  } Let us go back to the continuous time process $\set{\BP(t):t\geq 0}$ and suppose for ``large $t$'', we sample a vertex $V_t$ uniformly at random from $\BP(t)$. Let $\fB(V_t) \leq t$ denote the time this vertex was born into the system, so that by time $t$, if we consider the $\age(V_t) = t - \fB(V_t)$ of this vertex then for any $a\geq 0$,
\begin{equation}
    \pr(\age(V_t) \geq a|\BP(t)) = \pr(V_t \in \BP(t-a)|\BP(t))=\frac{|\BP(t-a)|}{|\BP(t)|} \to \exp(-a),
\end{equation}
where in the last assertion, we have used \eqref{it:exp}. Now since in continuous time all the dynamics is independent,  the following seems like a believable guess:

\noindent{\bf Guess:} The fringe distribution of the $\set{\cT_k:k\geq 1}$ should converge to the genealogical tree of of the branching process $\BP$ observed for a random exponential mean one amount of time. 

This is true, and in fact the above implies something much deeper. 

\begin{thm}[\cite{aldous-fringe,jagers-nerman-1,jagers-nerman-2,jagers-ctbp-book}]
    The sequence of random trees $\cT_n \probcrf \cT_\infty $ where $\cT_\infty$ is the random sin-tree whose fringe distribution $\varpi$ is constructed as follows:
    \begin{enumeratei}
        \item Generate $T\sim \exp(1)$.
        \item Generate a continuous time branching process $\BP$ with offpsring distribution given by a Poisson process of rate one, independent of $T$. 
    \end{enumeratei}
    Let $\BP(T)$ geneological tree of the branching process observed for $T$ units of time then  $\varpi$ is the distribution on $\bT$ of this random finite tree.
\end{thm}


Spurred by the above discussion, it turns out much of the above discussion can be extended to much more general tree models.   Consider the general random tree model sometimes called \textit{non-uniform random recursive trees} \cite{szymanski1987nonuniform}.
Fix $n \ge 1$ and an \textit{attachment function} $f:\{0,1,2\dots\} \to  (0, \infty)$.
{A sequence of random trees $\set{\cT_k: 1\leq k \leq n}$ is grown as follows ($\cT_k$ has $k$ vertices labelled by the integers $[k] := \{1, \dots, k\}$). For $k=1$, $\cT_1$ has one vertex, which we call the ``root.''}
For fixed $k \ge 2$, $\cT_k$  is constructed conditional on $\cT_{k-1}$ as follows.
{A new vertex, $k$, is born into the system and attaches to a previously existing vertex $v \in [k-1]$ with probability proportional to $f(\text{deg}(v))$, where $\text{deg}(v)$ denotes the number of children of $v$ (which is one less than its graph degree in $\cT_{k-1}$).}
Thus,
\begin{equation*}
{\mathbb{P}\left(k \text{ attaches to } v \in [k-1] \ \vert \ \cT_{k-1}\right) := \frac{f(\text{deg}(v))}{\sum_{u=1}^{k-1} f(\text{deg}(u))}.}
\end{equation*}
The vertex that $k$ selects is called the ``parent" and the edge is directed from the parent to the new ``child" vertex.
The case of $f(\cdot)\equiv 1$ corresponds to the above model we described and is now known known under the phrase {\textit{random recursive trees}} \cite{smythe1995survey}.
The specific case of \emph{linear preferential attachment} when $f$ is affine, namely $f(k) = k+1+\beta,~k\geq 0,$ for a fixed parameter $\beta$ was considered in \cite{barabasi1999emergence} to provide a generative story for heavy tailed degree distributions of real networks. 

Then it turns out, that each of these models can be embedded in continuous time in corresponding continuous time branching processes $\BP_{\text{\tt model}}$,  as we did for the random recursive tree model. Further,  under technical conditions on the attachment function $f$, each of the models have a (model dependent) constant $\lambda  = \lambda_{\text{\tt model}}>0$ (called the Malthusian rate of growth of this branching process) such that $|\BP_{\model}(t)| \sim W_{\text{\tt model}} e^{\lambda t}$. Thus the above heuristic once again leads to guesses as well as a proof path to show this class of random tree models converge to limiting (model dependent) random {\tt sin}-trees. Further discussion is postponed to Section \ref{sec:dynamic-res}.  

\subsubsection{Networks: \erdos and the Configuration model}\label{staticlwl}
While all the results here are well known, a description of the history and proofs can be found in \cite[Chapter 2]{van2023random}. 
Recall the intuition developed related to the \erdos random graph in Section \ref{sec:lwc-gfw}. The following should not come as a surprise.
\begin{thm}
    Fix $\lambda >0$ and consider the sequence of random graphs $\set{\cG_n:n\geq 1}$ with $\cG_n \sim \ERRG(n,\lambda/n)$. Let $\cG_\infty \sim \BP(\pois(\lambda))$ denote a Galton-Watson Branching process with $\pois(\lambda)$ offspring distribution. Then $\cG_n \probLWC \cG_\infty$ as $n\to\infty$.   
\end{thm}
This is also a good place to describe functionals that local weak convergence might {\bf not} be able to capture {\bf without more work or conditions.} Consider the sequence of network models $\set{\tilde{\cG}_n:n\geq 1}$ via:
\begin{enumeratea}
    \item $\tilde{\cG_n}$ has vertex set $[2n]$ partitioned as $[1:n]$ and $[n+1:2n]$. 
    \item Generate independent $\ERRG(n,4/n)$ on each of these two vertices.
    \item Select $U_n^{\sss(1)},U_n^{\sss(1)}$ uniformly at random from $[1:n]$ and $[n+1:2n]$ and connect them by a single edge. 
\end{enumeratea}
Then it is possible to check that $\tilde{\cG}_n \probLWC \cG_\infty$ as $n\to\infty$ as in the above result, however this sequence of networks are {\bf qualitatively} different from the $\ERRG$ model; for example with some strictly positive probability as $n\to\infty$, these networks have two disjoint giant components (components of the same order as the system size $2n$) as $n\to\infty$. This should not be surprising as local weak convergence is not built to reveal asymptotics for global functionals (such as the connectivity properties of the network). What will be surprising is, with a lot of hard work, one can push this notion to give detailed asymptotics for global functionals in a wide variety of models. 

The next example describes a mechanism for generating network models such that (in the limit), the empirical degree distribution converges to a pre-specified pmf $\vp:= \set{p_k:k\geq 0}$. There are various variants of the following model, known as the \emph{configuration model}, and we will pick one specific construction. For network with vertex set $[n]$: 
\begin{enumeratea}
    \item Generate $\vd_n:= \set{d_i: i\in [n]}$ iid with distribution $\vp$. Assume $\ell_n = \sum_{i\in [n]} d_i$ is even else add one to $d_n$. 
    \item Imagine every vertex $i \in [n]$ having $d_i$ half edges, where two half-edges have to be paired to create a full edge. Label  $\ell_n$ half edges in some arbitrary fashion from $1,2,\ldots, \ell_n$.
    \item Start half-edge labelled one and pair this one of the remaining $\ell_n-1$ half-edges uniformly at random to form a full edge and remove both of these half-edges; next pick the smallest indexed unpaired half-edge and pair this one of the remaining $\ell_n-3$ remaining half-edges to form the second full edge and remove both of these from the collection of unpaired half-edges; continue till all half-edges have been paired.  
\end{enumeratea}
Write $\cG_n \sim \CM_n(\vp)$ for the random network  obtained at the end of this process; in principle this will  be an \emph{multi-graph} as there might potentially be multiple edges and self-loops. 

\noindent {\bf Intuition for asymptotics:} In the large network $n\to\infty$ limit,  suppose we pick a vertex $V_n$ uniformly at random from this graph; then the number of neighbors (namely initial half-edges) of this vertex should have distribution $\vp$. Let us follow one of these half-edges of $V_n$ to understand the degree of the vertex whose half-edge connects to this half-edge completing a full edge; conceptually vertices with degree $k$ (of which there are a density $p_k)$ are $k$ times more likely to complete this edge than a vertex of degree $1$. Thus it is not hard to convince oneself that the probability that the degree of the vertex that completes this edge is $k$ should be proportional to $kp_k$. However after completing this edge, this new vertex has $k-1$ half-edges to connect to subsequent vertices. This motivates the following definitions. 

\begin{defn}[Unimodular branching processes with input pmf $\vp$]
\label{def:unimodular-bp}
    Assume the pmf $\vp$ has finite mean $\mu = \sum_{k=0}^\infty kp_k < \infty$. Define the size-biased pmf $\vp^\circ$ given by:
    \[p_k^\circ = \frac{(k+1)p_{k+1}}{\mu}, \qquad k\geq 0. \]
     Consider the branching process $\BP^{\circ}(\vp, \vp^\circ)$ started from a single root $\rho$ where the root has offspring distribution $\vp$ but every subsequent generation has offspring distribution $\vp^\circ$. 
\end{defn}

While the following is well known in the literature, a full proof can once again be found in \cite[Chapter 2]{van2023random}. Here Let $D\sim \vp$ will denote a generic random variable with distribution given by the degree pmf generating $\CM_n(\vp)$. 

\begin{thm}
\label{thm:config-res}
    Assume $\E(D^2) < \infty$. Then the sequence of models $\set{\cG_n:n\geq 1}$ with $\cG_n\sim \CM_n(\vp)$ satisfy $\cG_n \probLWC \BP^\circ(\vp, \vp^\circ)$ as $n\to\infty$ in the sense of Definition \ref{def:lwc-prob}. 
\end{thm}





\section{Power of local weak convergence}
\label{sec:power}

\subsection{An illustrative example: Spectral distribution of adjacency matrices for trees}

Recall that $\mu^{\sss(n)}_{\EMP}$ denoted the empirical distribution of the eigen-values of $\cT_n$. 

\begin{thm}[\cite{bhamidi2012spectra}]
Consider a sequence of trees converging in fringe sense to a random infinite {\tt sin}-tree. Then there exists a model dependent, non-random, probability measure $\mu_{{\model}}$ such that $\mu^{\sss(n)}_{\EMP} \probc \mu_{{\model}}$ as $n\to\infty$. 
\end{thm}

\noindent{\bf Rationale for why this might be true:} Recall that $\cM_{\pr}(\bR)$ denotes the space of probability measures on $(\bR, \cB(\bR))$. For $\mu \in \cM_{\pr}(\bR) $, write $\supp(\mu)$ for the support of $\mu$.   For complex $z\in \bC, z\notin \supp(\mu)$, define the Stieltjes transform at $z$,
\begin{equation}
\label{eq:stj-tr}
s_\mu(z):= \int_{\bR} \frac{1}{x-z} d\mu(z). 
\end{equation}
One can check that if $z = x+i y$ with $y>0$ then $|s_\mu(z)| < 1/y$. It seems believable that if one can show, that there is a (model dependent) deterministic function $s_{\model, \infty}(z)$ such that for $z\in \vA\subseteq \bC$,   $s_{\mu^{\sss(n)}_{\EMP}}(z) \probc s_{\model, \infty}(z)$, where the domain of convergence $\vA$ is ``large enough'', then this should at least augur convergence of the empirical distribution to a limit deterministic measure $\mu_{\model}$ whose Stieltjes transform is given by $s_{\model, \infty}(z)$. 

 Now for $\mu^{\sss(n)}_{\EMP}$, with $\mbox{Tr}$ denoting the trace of a matrix, 
\begin{equation}
\label{eqn:stiel-def}
    s_{\mu^{\sss(n)}_{\EMP}}(z) = \frac{1}{n}\mbox{Tr}(\vA_n-zI)^{-1} := \frac{1}{n}\sum_{v\in \cT_n} R_{(vv), \cT_n}(z),
\end{equation}
where $R_{(vv), \cT_n}(\cdot)$ is sometimes referred to as the \emph{resolvent} of vertex $v$ in the tree $\cT_n$. 
Let $\cN(v)$ denote the number of children of $v$ in $\cT_n^-(v)$ and denote these as $\set{c_i(v):1\leq i\leq \cN(v)}$ and the corresponding subtrees hanging from these vertices by $\cT_{n,i}^-(v)$ (viewed as a rooted tree at $c_i(v)$). Let $\rho(v)$ denote the parent of $v$ (namely the unique node on the path to the root of $\cT_n$) and let $\cT_n^{+}(v)$ denote the tree rooted at $\rho(v)$ obtained by deleting the edge $(v, \rho(v))$ and the fringe tree below $v$.

    In the setting of trees, it turns out via Schur-decomposition that the resolvent at $v$ in the original tree $\cT_n$ can be decomposed as, 
\begin{equation}
 \label{eqn:res-decom}
 R_{(vv), \cT_n}(z) = \frac{1}{-z+ \sum_{i=1}^{\cN(v)} R_{(c_i(v) c_i(v)), \cT_{n,i}^-(v)}(z) + R_{(\rho(v), \rho(v)), \cT_n^{+}(v)}(z)}.
 \end{equation} 
Fix $z\in \bC$ with imaginary part $Im(z)> 1$ and iterate the above expansion $d$ times, resulting in a continued fraction type expansion, where we have $d$ terms dependent on the fringe expansion upto distance $d$ and a remainder term that depends on the graph beyond distance $d$ from $v$.  Intuitively it seems clear that, if one has extended fringe convergence as in Def. \ref{eqn:deg-convg-fr}, then when $Im(z) >1$, this should thus imply that $s_{\mu^{\sss(n)}_{\EMP}}(z) \probc s_{\model, \infty}(z)$. Converting this idea into a rigorous proof takes a lot of work and can be found in \cite{bhamidi2012spectra}. 

\noindent {\bf What did we learn?} Conceptually the spectral distribution of the adjacency matrix is a {\bf global} functional and depends on the entire network. However if one has local weak convergence, which a priori only gives information on asymptotics of local functionals, one can exploit it to also establish convergence of global functionals, provided that the effect of `far away' vertices on the functional value of a given vertex diminishes in a quantifiable way.

\subsection{Random matrices for general random graph models}
Now let us consider the case of general (not-necessarily tree) networks that still converge in the local weak sense to trees. Our goal is to describe a small sub-collection of the foundational results in \cite{bordenave2010resolvent}.  We will mainly focus on the configuration model driven by pmf $\vp$  (owing to the assumption in Theorem \ref{thm:config-res} we assume $\vp$ has finite second moment). We let $\BP(\vp^\circ)$ denote a (standard) branching process with offspring distribution $\vp^\circ$ while $\BP^\circ(\vp,\vp^\circ)$ denotes the unimodular branching process as in Def. \ref{def:unimodular-bp}; similarly define $D\sim\vp, D^\circ\sim \vp^\circ$. Let $\cG_n\sim \CM_n(\vp)$ and let $\vA_n, \mu^{\sss(n)}_{\EMP}$ as before denote the corresponding adjacency matrices and empirical spectral distribution respectively.  The Stieltjes transform as in \eqref{eqn:stiel-def} will be our main tool again so let us give a name to the class of such functions, with $\bC_+:=\set{z \in \bC, \fI(z) >0}$:
\begin{equation}
    \label{eqn:h-def}
    \cH = \set{f: \bC_+ \to \bC_+, f \text{ holomorphic }, |f(z)|\leq \frac{1}{\fI(z)}}.
\end{equation}
Paraphrasing a small sub-class of many amazing results:
\begin{thm}[{\cite{bordenave2010resolvent}}]
\label{thm:bord-lel}
    \begin{enumeratea}
        \item Let $\cP(\cH)$ denote the class of probability measures on $\cH$. There is a unique probability measure $\bQ \in \cP(\cH)$ such that $Y\sim \bQ$ satisfies the following recursive distributional equation:
        \[Y(z) \stackrel{d}{=} -\frac{1}{z+ \sum_{i=1}^{D^\circ} Y_i(z)}\]
where $\vY:=\set{Y_i(\cdot):i\geq 1}$ are \emph{iid} with distribution $\bQ$ and independent of $D^\circ \sim \vp^\circ$. 
\item Let $D\sim \vp$ independent of $\vY$ and define the random element $X\in \cH$ via,
\[X(z) = -\frac{1}{z+ \sum_{i=1}^D Y_i(z)},\]
Then $\mu^{\sss(n)}_{\EMP} \probc \mu_{\infty,\vp}$ where $\mu_{\infty,\vp}$ is the unique probability measure on $\cP(\bR)$ with Stieltjes transform $s_{\infty, \vp}(z) = \E(X(z))$. 
    \end{enumeratea}
\end{thm}
\noindent{\bf Illustrative example:} Fix $k\geq 3$ and consider the degree distribution $\vp \sim \delta_k$, this then results in the random $k$-regular graph, which converges in the LWC sense to a rooted tree where the root has degree $k$ and all the remaining individuals have degree $k-1$. Suppose, as a guess, one searches for deterministic solutions $Y(\cdot)$ in (a) above. This results in $Y(z) = -(z+(k-1)Y(z))^{-1}$. Solving, one finds that this is the Stieltjes transform of the semi-circular law with radius $2\sqrt{k-1}$. Using this in (b) shows that the empirical spectral distribution of random $k$-regular graph converges to the probability measure with density 
\[f(x) = \frac{k}{2\pi} \frac{\sqrt{4(k-1) -x^2}}{k^2-x^2}, \qquad x\in [-2\sqrt{k-1}, 2\sqrt{k-1}]. \]
 
\begin{figure}
    \centering
    \includegraphics[scale=.25]{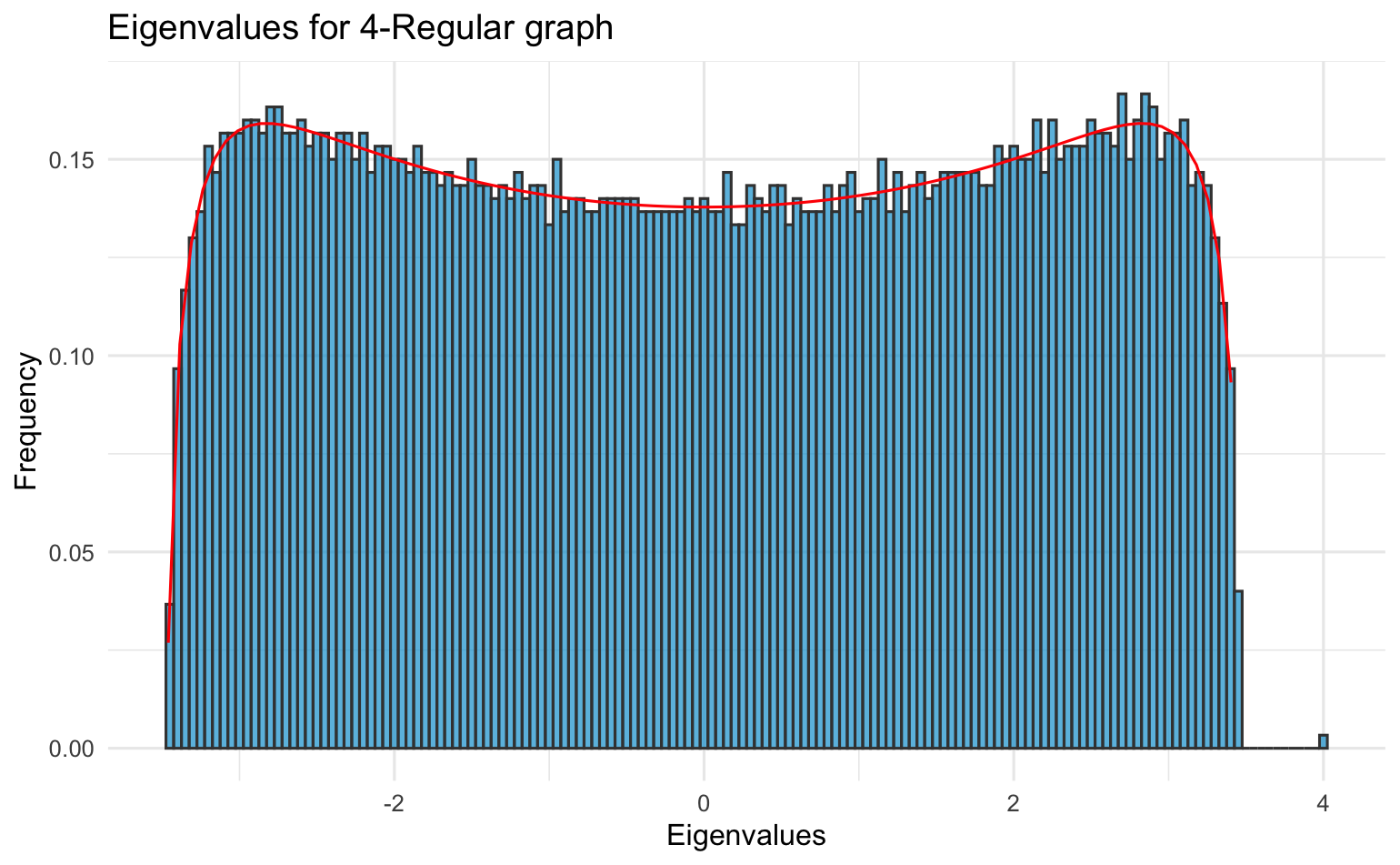}
    \caption{A simulation with $k=4$ and network size $n=6000$ showing the empirical spectral distribution with the theoretical limit.}
    \label{fig:esd-empirical}
\end{figure}

\noindent{\bf Rationale for why this might be true:} There are four steps in the proof in \cite{bordenave2010resolvent}:
\begin{enumeratei}
    \item {\bf Local weak convergence implies convergence:} Using general functional analytic techniques, the authors are able to show that local weak convergence (not even necessarily for trees) with some minor technical conditions automatically implies convergence of the Stieltjes transform with the limit Stieltjes transform $s_\infty(z) = \lim_{n\to \infty} \E(R_{V_n, V_n}(z))$, where $R_{V_n, V_n}(\cdot)$ denotes the resolvent of a uniformly selected vertex $V_n$. 
    \item {\bf Specific analysis when one has convergence to $\BP^{\circ}$:} Now let us consider the case of the configuration model where one has convergence to $\BP^\circ$. Let us start with an additional simplification: let us try to understand what happens to the adjacency matrix of the standard branching process $\BP_m(\vp^\circ)$ grown to $m$ levels say, with adjacency matrix $\vA_m$, in particular what happens to the resolvent $Y_\rho^{\sss m}(z):= R_{\rho, \rho}(z) = \tr(\vA_m - z\vI)^{-1}_{\rho, \rho}$. Let $D_\rho^{\circ}\sim \vp$ denote the number of children of the root and write these as $\set{v_1,v_2, \ldots, v_{D_\rho^\circ}}$. Note that the subtrees below each of these vertices constitute independent Branching processes run upto $m-1$ generations, see Figure \ref{Tree}. 
    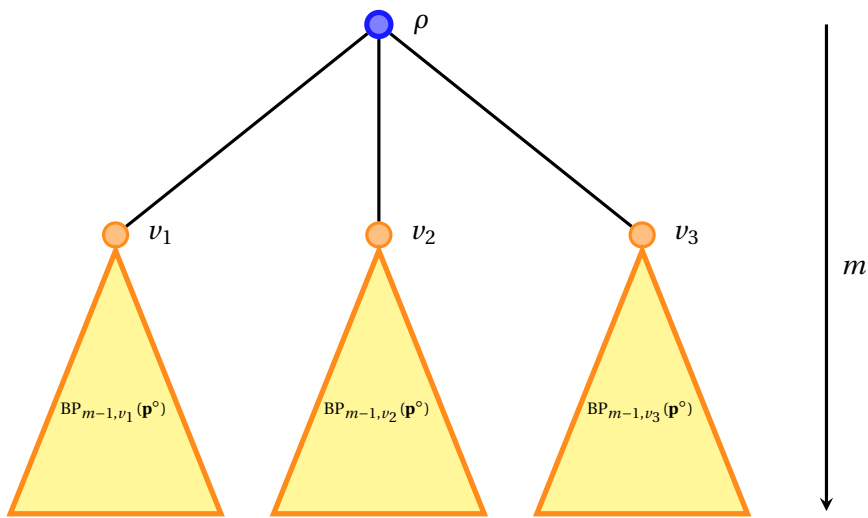
\begin{figure}
    \centering
    \begin{tikzpicture}[scale=.35, auto=left, every node/.style={circle, minimum size = 0.3 cm}]
        \node[draw = orange!90, fill = orange!50,very thick, label=right:{$v_2$}] (n3) at (10,10) {};
        \node[draw = orange!90, fill = orange!50,very thick, label=right:{$v_3$}]  (n4) at (20,10) {};
        \node[draw = orange!90, fill = orange!50,very thick, label=right:{$v_1$}]  (n2) at (0,10) {};
        \node[draw = blue!90, fill = blue!50, label=right:{$\rho$}]  (n1) at (10,18) {};
        
         \draw[draw = orange!90, fill = yellow!50] (n2)++(0,-0.6) -- ++(-4,-10) -- ++(8,0) -- node[pos = 0.57] {{\tiny $\BP_{m-1,v_1}(\vp^\circ)$}} cycle;
        \draw[draw = orange!90, fill = yellow!50] (n3)++(0,-0.6) -- ++(-4,-10) -- ++(8,0) -- node[pos = 0.57] {{\tiny $\BP_{m-1,v_2}(\vp^\circ)$}} cycle;
        \draw[draw = orange!90, fill = yellow!50] (n4)++(0,-0.6) -- ++(-4,-10) -- ++(8,0) -- node[pos = 0.57] {{\tiny $\BP_{m-1,v_3}(\vp^\circ)$}} cycle;
        \draw[black,very thick](n1) -- (n2) -- cycle;
        \draw[black,very thick](n1) -- (n3) -- cycle;
        \draw[black,very thick](n1) -- (n4) -- cycle;
        \draw[->, very thick, >=stealth] (n1)++(17,0) -- ++(0,-18.5) node[midway,right] {$m$};
    \end{tikzpicture}
    \caption{Figure of $\BP^\circ(\vp^\circ)$ with root degree $D^\circ_\rho = 3$}
    \label{Tree}
\end{figure}

    Let $\vA_{m-1}^{v_i}$ denote the adjacency matrices of each of these subtrees and let 
    \[Y_{v_i}^{m-1}(z) = \tr(\vA_{m-1}^{v_i} - z\vI)^{-1}_{v_i, v_i} \]
    An application of Schur-decomposition as in the previous Section results in the equation:
    \begin{equation}
        \label{eqn:443}
        Y_\rho^{\sss m}(z) = - \frac{1}{z+ \sum_{i=1}^{D_\rho^{\circ}}  Y_{v_i}^{m-1}(z) }.
    \end{equation} 
    \item {\bf Recursive distributional equations and uniqueness:}
    Consider the map $\Psi: \cP(\cH) \to \cP(\cH)$ defined as follows: For $\bQ\in \cP(\cH)$ define $\Psi(\bQ)$ to be the distribution of the random variable $\tilde{Z}$ obtained as follows:  
    \begin{enumeratea}
        \item Generate $\vZ = \set{Z_i(\cdot):i\geq 1}$ iid with distribution $\bQ$ and independent of $D^\circ\sim \vp^\circ$.
        \item Let
        $\tilde{Z}(z) = -({z+ \sum_{i=1}^{D_\rho^{\circ}}  Z_i(z) } )^{-1}$. 
    \end{enumeratea}
    Say that $\bQ$ is a solution of the recursive distribution equation if $\Psi(\bQ) = \bQ$. Suppose we were able to show that this recursive distributional equation has a {\bf unique} solution. Then  \eqref{eqn:443} should (after some work) imply that as $m\to\infty$, $Y^m(\cdot)$ converges to this unique solution as $m\to\infty$. The authors in \cite{bordenave2010resolvent} are able to prove this uniqueness via showing that this mapping is a contraction and then using Banach fixed point theorem. 
    \item {\bf Relating this to $\CM_n$:} Redoing the Schur decomposition analysis but now 
   to understand $R_{V_n, V_n}(\cdot)$ for a randomly selected vertex in $\CM_n$, and connecting it to the above analysis of $\BP^\circ(\vp, \vp^\circ)$ via local weak convergence as described before, it should now be believable that $R_{V_n, V_n}(z) \approx X(z)$ for large $n$ where $X(z)$ satisfies the distributional equation in (b) of the Theorem. 
\end{enumeratei}

\noindent {\bf What did we learn?} Here the local computations not just proved the existence of the limit (of the empirical spectral distribution) but owing to the recursive construction of the limit using independent and identically distributed components (here, independent branching processes hanging off the root), this resulted in explicit information about the limit via \emph{computations} carried out on the limit object. This phenomenon related to local weak convergence was termed the \emph{objective method} in \cite{aldous-steele-obj} where local weak convergence allows one to:
\begin{enumeratea}
    \item Prove asymptotics for functionals of interest in the large network limit;
    \item Obtain information on the limit via recursive computations performed on the limit object. 
\end{enumeratea}

\subsection{Interacting particle systems on random graphs}
Recall the motivation in Section \ref{sec:gibbs-distr-intro}. Once again our goal will be to understand properties of the Gibbs distribution on a sequence of configuration model driven graphs $\cG_n \sim \CM_n(\vp)$ as $n\to\infty$. Let $\mu_{\cG_n}(\cdot)$ denote the probability measure in \eqref{eqn:gibbs}. Recall that in the previous sections, the first step was to relate such questions to a functional (the Stieltjes transform for example) that we could use to pave the way for understanding asymptotics. In this setting, the normalizing factor $Z_n(\beta,B)$, which is a function of the inverse temperature $\beta$ and external field $B$ seems like a natural candidate. Since the number of potential configurations (terms in the summand) is exponential, the following object:
\begin{equation}
    \phi_n(\beta, B) = \frac{1}{n} \log{Z_n(\beta, B)},
\end{equation}
seems like a natural starting point. Further note that appropriate derivatives of this function give us important information, for example, the mean magnetization is given by, 
\begin{equation}
    \label{eqn:133}
    \frac{\partial \phi_n(\beta, B) }{\partial B} = \sum_{\vx\in \set{-1,1}^n} \frac{\sum_{i\in [n]} x_i}{n}\mu_{\cG_n}(\vx) = m_{\cG_n}(\beta, B).
\end{equation} 
Thus we try to understanding asymptotics of this object as $n\to\infty$. Our initial goal in this Section is to describe the foundational results in \cite{dembo2010ising} describing these asymptotics in the \emph{ferromagnetic} ($\beta \geq 0$) regime. Since the model is ``symmetric'' with respect to $B$ (in the sense that $\phi_n(\beta, B) = \phi_n(\beta, -B)$), throughout this Section we will assume that $B>0$.   We will now describe the foundational result in \cite{dembo2010ising}, the reader should contrast the result and the preceding setup with Theorem \ref{thm:bord-lel}. First we need some notation. Throughout this Section we assume $\beta\geq 0$. Let $\tanh,\atanh$ denote the tan hyperbolic and inverse tan hyperbolic functions respectively. Define the function, 
\begin{equation}
    \label{eqn:236}
    \fE(\beta, h):= \atanh[\tanh(\beta)\tanh(h)]
\end{equation}
Consider the functional $\Psi_{\beta, B}: \cP(\bR) \to \cP(\bR)$ from the space of probability measures on $\bR_+$ to itself given as the distribution of the random variable $X$ obtained via:
\begin{enumeratea}
    \item Generate $\set{Y_i:i\geq 1}\sim_{iid} \pr$ independent of $D^{\circ}\sim \vp^\circ$. 
    \item Let $X = B+ \sum_{i=1}^{D^\circ -1} \fE(\beta, Y_i)$. 
\end{enumeratea}
\begin{thm}{\cite{dembo2010ising}}
    \label{thm:dem-mont} Consider the Ising model as in \eqref{eqn:gibbs} on $\CM_n(\vp)$ with finite second moments and with $\beta \geq 0, B\geq 0$. Then, 
    \begin{enumeratea}
        \item For $\beta\geq 0$ the operation $\Psi$ is monotone increasing (under the stochastic ordering operation) and there is a unique fixed point $\tilde{\pr}_{\beta, B}$ supported on $[0,\infty)$ satisfying $\tilde{\pr}_{\beta, B} = \Psi(\tilde{\pr}_{\beta, B})$ i.e. for $Y,\vY:=\set{Y_i:i\geq 1}\sim_{iid} \tilde{\pr}_{\beta, B}$ independent of $D^\circ$, 
        \[Y \stackrel{d}{=} B+ \sum_{i=1}^{D^\circ -1} \fE(\beta, Y_i).\]
        Call this the distribution of ``local fields''.
        \item Let $D\sim \vp$ independent of $\vY$. Then $\phi_n(\beta, B) \probc \phi_\infty(\beta, B)$ where
      \begin{align*}
            \label{eqn:249}
            \phi_\infty(\beta,B) = &\frac{\mu}{2} \log\cosh(\beta) -\frac{\mu}{2} \E[\log[1+\tanh(\beta)\tanh(Y_1)\tanh(Y_2)]] \\
            &+ \E\log\bigg\{e^B \prod_{i=1}^D[1+\tanh(\beta)\tanh(Y_i)]+e^{-B} \prod_{i=1}^D[1-\tanh(\beta)\tanh(Y_i)]\bigg\}.
            \end{align*}
       
    \end{enumeratea}
\end{thm}

\noindent{\bf Rationale for why this might be true:} The proof of this result is significantly more complicated so we will only be able to describe a small aspect of the proof. We will write $\inpr{\nu}{f} = \sum_x f(x)\nu(x)$ for the expectation operation of a function with respect to a measure on a finite set. There are four major ideas in the proof: 

\begin{enumeratea}
\item {\bf Monotonicity in the Ferromagnetic regime:} In the ferromagnetic regime, it turns out that there are fundamental monotonicity properties of the model that one can take advantage. Two important threads in this direction are: 
\begin{enumeratei}
    \item {\bf Monotonicty with respect to local magnetic fields:} For a graph $\cG$ and vertex dependent non-negative constants $\vB = \set{B_u:u\in V(\cG)}$ (called local field) consider a generalized version of the model in \eqref{eqn:gibbs} on a graph via:
    \begin{equation}
\label{eqn:gibbs-ext}
\mu_{\vB, \beta}(\vx) = \frac{1}{Z(\beta, \vB)} \exp\left(\beta \sum_{(u,v)\in E(\cG)} x_u x_v + \sum_{u\in V} B_u x_u\right)
\end{equation}
For two local field parameters $\vB, \vB^\prime$ with $0\leq B_u\leq B_u^\prime$ for all $u\in \cG$ and inverse temperature parameters $0\leq \beta\leq \beta^\prime$  Griffith's inequality says that for any subset $U\subset V(\cG)$, 
$0\leq \inpr{\mu_{\beta, \vB}}{\prod_{i\in U} x_i} \leq \inpr{\mu_{\beta^\prime, \vB^\prime}}{\prod_{i\in U} x_i} $. 
\item {\bf Monotonicity with respect to boundary conditions: } To reduce notational overhead, we will suppress dependence on $\beta, \vB$ of the measures in question unless required.  Fix a connected graph $\cG$ and a strict subset of vertices $U\subset \cG$ and let the boundary $\partial U$ denote the collection of vertices in $U$ that have connections to at least one vertex in $\cG \setminus U$. Suppose one was interested in the marginal distribution of the Ising model on the subgraph spanned by vertices in $U$. Then intuitively there are two Ising type models that ``lower'' and ``upper'' bound this distribution, the Ising model on $U$ with free ($\mu^{U, 0}$) and all $+$ ($\mu^{U, +}$) boundary conditions respectively: 
\begin{align*}
    \mu^{\sss U, 0}(\vx_U) &= \frac{1}{Z^{\sss U,0}(\beta,\vB)} \exp\bigg(\beta \sum_{(u,v)\in E(U)}x_u x_v + \sum_{u\in U } B_u x_u\bigg),\\
    \mu^{\sss U, +}(\vx_U) &= \frac{1}{Z^{\sss U,+}(\beta,\vB)} \exp\bigg(\beta \sum_{(u,v)\in E(U)}x_u x_v + \sum_{u\in U } B_u x_u\bigg)\ind\set{\vx_{\partial U} \equiv \boldsymbol{+1}},
\end{align*}
i.e. when attempting to compute the marginal distribution of the Ising model on $U$, the first measure assumes no influence from the outside (in the ferromagnetic regime one would imagine some sort of ``positive influence'' from the boundary) while the second measure assumes maximal positive influence from the boundary in the sense that all these vertices are set to have spins $+1$. 
\end{enumeratei}
    \item {\bf Connecting this to local weak convergence on trees and vanishing impact of boundary conditions:} While we are interested in $\phi_n$, \eqref{eqn:133} already shows how \emph{derivatives} of this object can be related to averages. After significant amount of work, what is shown in \cite{dembo2010ising} is that to prove the main result, it is enough to show that for each $\beta \geq 0$,
    \[\frac{\partial \phi_n(\beta, B)}{\partial \beta} \probc \frac{\partial \phi_\infty(\beta, B)}{\partial \beta}. \]
    Writing $\mu_{\cG_n}$ for the Ising model on the graph $\cG_n$, direct computation shows that, 
    \[\frac{\partial \phi_n(\beta, B)}{\partial \beta} = \frac{1}{n} \sum_{v\in [n]} \sum_{u\sim v} \inpr{\mu_{\cG_n}}{x_u x_v} = \E_n(\sum_{u \in \cN( V_n)} \inpr{\mu_{\cG_n}}{x_u x_{V_n}} ),\]
    where $\E_n$ denotes expectation conditional on the graph $\cG_n$ and thus the only randomization is over the choice of the random vertex $V_n$ sampled uniformly from the graph. Now fix any $t\geq 1$ and as before let $B_{\cG_n}(V_n, t)$ denote the ball of radius $t$ around $V_n$. Conditional on $V_n$, let $\mu^{\sss B_{\cG_n}(V_n, t),0}, \mu^{\sss B_{\cG_n}(V_n, t),+} $ denote the Ising models with free and + boundary conditions on this neighborhood and let $\E_n^{0,t}, \E_n^{+,t}$ denote expectations of functionals over choices of $V_n$. Then the monotoncity described in the previous step implies for each $t$, 
    \[\E_n^{0,t}\bigg(\sum_{u \in \cN( V_n)} \inpr{\mu^{\sss B_{\cG_n}(V_n, t),0}}{x_u x_{V_n}} \bigg)\leq \frac{\partial \phi_n(\beta, B)}{\partial \beta} \leq \E_n^{+,t}\bigg(\sum_{u \in \cN( V_n)} \inpr{\mu^{\sss B_{\cG_n}(V_n, t),+}}{x_u x_{V_n}} \bigg).\]
    Now local weak convergence implies that as $n\to\infty$, the left hand side and right side converge to the following: Consider $\BP(\rho,t):= \BP^\circ(\vp, \vp^\circ)$ grown upto generation $t$ from root $\rho$. Consider the Ising model on this with the same parameters with free and + boundaries and let $\E_\infty^{0,t}, \E_\infty^{+,t}$ denote expectations over the branching process measure. Then the left hand side should converge to
    \[\E_n^{0,t}\bigg(\sum_{u \in \cN( V_n)} \inpr{\mu^{\sss B_{\cG_n}(V_n, t),0}}{x_u x_{V_n}} \bigg) \probc \E_\infty^{0,t}\bigg(\sum_{u \in \cN(\rho)} \inpr{\mu^{\sss \BP(\rho, t),0}}{x_u x_\rho} \bigg)\]
and similar expressions for the right hand side. Further it is shown that as $t\to\infty$, the effect of the boundary conditions on the above expectations vanish (a highly consequential finding). An example of quantification of such a result is as follows: consider the branching process $\BP^\circ(\rho,m)$ as above upto generation $m$ and the corresponding Ising model $\mu^{m}$ say with corresponding free and + measures $\mu^{m, 0}, \mu^{m,+}$ .  Suppose for some $t< m$ we are interested in the marginal distribution $\mu_U^{m,0}, \mu_U^{m,+}$  of the Ising model on a subset of vertices $U\subset \BP^\circ(\rho,t)$, in particular the impact of the above boundary conditions. Then \cite[Theorem 4.2]{dembo2010ising} shows that (under the assumptions of the Theorem) there are constants $M, C$ such that the expected total variation distance satisfies,
\[\E(||\mu_U^{m,0} - \mu_U^{m,+}||_{TV}) \leq \frac{\E(C^{|\BP(t)|})}{m-t}.\]

    \item {\bf Pruning argument for the Ising model on trees and the origin of local fields:}
    This sequence of arguments leads to trying to understand the Ising model on $\BP^{\circ}(\rho,m)$ for ``large $m$'' (where the boundary conditions are inconsequential) and in particular asymptotics of expectations of the sum taken over products of root vertex spin with its neighbors. Now we are back at Figure \ref{Tree}, the Ising model on this tree, and in particular the \emph{marginal distribution of the root and its neighbors}. If we had a tractable description for this, we are in business. This is where another fundamental finding related to Ising models for trees in \cite{dembo2010ising} comes to the rescue.
    Their general result says that the marginal distribution of the root and its neighbors can be \emph{recursively} constructed as follows:
    \begin{enumeratei}
        \item For each neighbor $v_i$ of the root,  consider the Ising model $\mu^{\BP_{m-1, v_i}}$ on the corresponding subtree $\BP_{m-1, v_i}(\vp^\circ)$ (this includes the vertex $v_i$). Let $\inpr{\mu^{\BP_{m-1, v_i}}}{x_{v_i}}$ denote the mean magnetization of the vertex $v_i$ in this model. 
        \item Then the marginal distribution on the root and its neighbors is again an Ising model with the same inverse temperature parameter $\beta$ and where the root has local field $B$ while each of its neighbors has local fields $h_i = \atanh(\inpr{\mu^{\BP_{m-1, v_i}}}{x_{v_i}}) $. 
        
    \end{enumeratei}
    The distribution of these local field random variables can be recursively computed and this is what leads to the main theorem. 
\end{enumeratea}
\begin{wrapfigure}{l}{6cm}
        \begin{tikzpicture}[scale=.22, auto=left, every node/.style={circle, minimum size = 0.3 cm}]
        \node[draw = orange!90, fill = orange!50,very thick, label=right:{{\small $v_2, B_2 = h_2$}}] (n3) at (10,10) {};
        \node[draw = orange!90, fill = orange!50,very thick, label=right:{{\small $v_3, B_3 = h_3$}}]  (n4) at (20,10) {};
        \node[draw = orange!90, fill = orange!50,very thick, label=right:{{\small $v_1, B_1 = h_1$}}]  (n2) at (0,10) {};
        \node[draw = blue!90, fill = blue!50, label=right:{$\rho$}]  (n1) at (10,18) {};
         \draw[black,very thick](n1) -- (n2) -- cycle;
        \draw[black,very thick](n1) -- (n3) -- cycle;
        \draw[black,very thick](n1) -- (n4) -- cycle;
        \end{tikzpicture}
\end{wrapfigure}

\noindent {\bf What did we learn?} This application in particular shows that it sometimes takes a lot of work to connect local asymptotics of models to asymptotics of global functionals, in particular, to quantify the impact of vertices ``far'' from a typical vertex which do not have significant impact in terms of the functional being considered. Further we arrived at these ``local fields'' $\set{h_i: 1\leq i\leq D^{\circ}_\rho}$ which quantified the impact of the (in the large $m\to\infty$ limit, infinite) subtrees below these children on the distribution of the root. In the next application, we will see similar local fields arising out of, at first sight, a calculation that does not make any sense!



\subsection{Random assignment problem}
Recall the motivation and setup of the random assignment problem in Section \ref{sec:prob-opt-mot}. Under moment conditions (and assumption of continuity)  on the edge cost distributions, it turns out that the limit depends only on the nature of the density at zero, so we will assume that the edge costs are \emph{iid} exponential random variables; further it will be convenient to take these to have mean $n$ (or rate $1/n$) variables. Before the nail in the coffin on this problem by Aldous in \cite{aldous1991asymptotic,aldous2001zeta}, \cite{goemans1993lower} had found a lower bound of 1.51 on the limit constant while \cite{coppersmith1999constructive} had found an upper bound of 1.94.

\begin{thm}[\cite{aldous2001zeta}]
    Consider the random assignment problem where edge costs are exponential mean $n$ random variables. Let $A_n$ denote the cost of the minimum assignment problem. Then:
    \begin{enumeratea}
        \item Let $\Xi = (\xi_1, \xi_2, \ldots,)$ be the points of a rate one Poisson point process and consider the recursive distributional equation  given by 
        \[X \stackrel{d}{=} \min_{1\leq i < \infty} (\xi_i - X_i),\]
        where $X, \set{X_i:i\geq 1}$ are \emph{i.i.d.} with the same distribution and independent of $\Xi$. Then the only solution to this distributional equation is the logistic distribution $X\sim f_X(\cdot)$ where the density is given by,
        \[f_X(x) = \frac{1}{(e^{x/2} + e^{-x/2})^2}, \qquad -\infty < x <\infty.\]
        \item Let $X_1, X_2\sim f_X(\cdot)$ and independent. Then 
        \[\frac{1}{n} \E(A_n) \to \int_0^{\infty} x\pr(X_1+ X_2 > x) dx = \frac{\pi^2}{6}.\]
    \end{enumeratea}
\end{thm}

\noindent{\bf Rationale for why this might be true:}

\begin{enumeratea}
    \item {\bf Asymptotics around a random job:} Writing $\pi_{opt}^{\sss(n)}$ for the optimal assignment, and switching from $c_{i,j}$ to $c(i,j)$ to easy the load on subscripts,  in the finite $n$ problem, note that by symmetry we see that  
    \[\frac{1}{n} \E(A_n) = \E\bigg(c({V_n,\pi_{opt}^{\sss(n)}(V_n) })\bigg).\]
    Thus once again this problem has been transformed into something related to what happens for a typical vertex, in this case the cost of the edge that vertex $V_n$ is matched to in the optimal assignment.
    \item {\bf Local weak convergence for Geometric graphs:} Let us think of the assignment problem as jobs being assigned to machines to facilitate explanation and by symmetry let $V_n = 1$. Recall that $\Xi = (\xi_1, \xi_2, \ldots)$ denotes a Poisson point process as above. Let $\set{Y_i:i\geq 1}$ be an iid sequence of exponential mean one random variables.   Now the least cost amongst possible machines to be assigned to this job is $\xi_{11}^{\sss(n)} := \min_{j\in [n]} c_{1j} \stackrel{d}{=} \xi_1 $ using properties of the exponential distribution. Using lack of memory property of the exponential distribution,  the second smallest cost
     has distribution \[\xi_{12}^{\sss(n)} \stackrel{d}{=}\xi_1 + \frac{n}{n-1} Y_2\approx \xi_2.  \]
    Similarly the 3rd least cost of possible machine assignments to job 1 has distribution,
    \[\xi_{12}^{\sss(n)} \stackrel{d}{=} \xi_{12}^{\sss(n)}  + \frac{n}{n-2} Y_3\approx \xi_3.\]
    Rigorous justification of the above facts proceeds via the R\'enyi representation of exponential order statistics.
    Thus, in the large $n$-limit, if we think of costs as edge lengths and explore the graph from job 1, then distances to machines are approximately distributed as a rate one Poisson point process. Now from each of these machines if we continue to explore the jobs closest to these machines (sequentially ``unfolding'' the geometric structure of costs from a single job), once again in the large $n$-limit they have Poisson rate one edge lengths, independent across machines. This suggests the following natural limit object. 
    \begin{defn}[Poisson weighted infinite tree (PWIT) \cite{aldous1991asymptotic}]
        Consider the infinite locally finite random tree constructed as follows: start with a single root $\rho$ and attach infinite number of children of this root where the children are assigned lengths according to a rate one Poisson process. Call this generation one.  Now recursively repeat the process, namely assuming we have constructed the process upto generation $k$,  for each vertex $v$ in any generation $k$ independently and identically repeat this process namely each of these individuals gives birth to individuals in generation $k+1$ with edge lengths according to a rate one Poisson process $\Xi_v$, independent across individuals. Call the random rooted tree with edge lengths the Poisson weighted infinite Tree (PWIT) and denote this by $\cT_\infty$.  Let $\vW$ denote the collection of edge lengths and $\vW(e)$ for the weight of a specific edge. 
    \end{defn}
    \item {\bf Assignments as matchings and the first approach that fails:} Recall how we arrived at the PWIT, this was via unfolding the cost structure around a typical job in the finite $n$ problem where we go from jobs in one generation to machines in the next generation and then to jobs in the generation after etc. Now any specific assignment of jobs to machines can, at least approximately (for large $n$) be thought of as a ``matching'' namely a collection of edges such that every vertex is part of exactly one edge. We can also define matchings $\cM$ on the PWIT, namely a collection of edges in $\cT_\infty$ such that every vertex is incident on exactly one edge.   Further, whilst at the far end of one's imagination, one can at least hope that something like the following is true:
\emph{
\begin{quote}
  There is an ``optimal'' matching on the PWIT such that as $n\to\infty$, $\E(c_n(1, \pi_{opt}^{\sss(n)}(1)))$ converges to $\E(\vW(\rho, \cM_{opt}(\rho)))$ namely the expected weight of the edge containing the root. 
\end{quote}
    }
    Now suppose we try the natural first approach to construct a matching on the PWIT which is the greedy matching approach $\cM_{\text{\tt greedy}}$: match $\rho$ with the vertex that is closest to it (which is at a distance of $\exp(1)$); then for all other children of $\rho$ match them with their closest child etc. Then it is obvious that under this scheme $\E(\vW(\rho, \cM_{\text{\tt greedy}})) =1$, however as mentioned above,  the lower bound for the limit constant before the result of Aldous was 1.51. Thus this cannot arise from the optimal matching. Conceptually it turns out that ``good'' matchings have to satisfy a spatial invariance property, in the sense that, going back to the definition of the PWIT, if we moved the root of the PWIT from the original root to one of its neighbors and then relabelled vertices, the matching should in some sense be invariant. It turns out that the greedy matching does not satisfy this property. In proving the existence of the limit constant, Aldous in \cite{aldous1991asymptotic} showed that $n^{-1}\E(A_n)$ converges to the the $\inf\E(\vW(\rho, \cM(\rho))$, where the infimum is taken over all spatially invariant matchings on the PWIT. 
    \item {\bf $\infty-\infty =$ Magic.  Origin of the Logistic distribution:} While trying to give an outline of the full and technical proof in \cite{aldous2001zeta} is beyond the scope of this article, let us describe one ``magical'' explanation (see \cite[Sec 4.2]{aldous2001zeta}) that shows the origin of the Logistic distribution above. Label the children of vertex $\rho$ from $\set{1,2,\ldots}$ and let $\cT_{j,\infty}$ denote the subtree rooted at $j$ which by construction has the same distribution as the PWIT. Let $C(\cT_\infty)$ denote the weight of the optimal matching on the PWIT (which is obviously $\infty$ but let us keep going). Let (this is \cite[Eqn 14]{aldous2001zeta}), 
    \begin{align*}
        X_\rho = &\text{cost of optimal matching on } \cT_\infty 
       -\text{cost of optimal matching on } \cT_\infty\setminus \set{\rho},
    \end{align*}
    and define similar variables for each of the children of $\rho$ namely $\set{X_j:j\geq 1}$ are the corresponding random variables for the subtrees $\cT_{j,\infty}$. By construction, if these random variables make sense in any universe, then they have the same distribution. For simplicity, on a set $\cU$ write $C(\cU)$ to denote the cost of the optimal matching on this set. Then note that 
    \[X_\rho = \min_{j\geq 1 } \bigg[W(\rho,j) + C(\cT_{j,\infty}\setminus \set{j}) + \sum_{i\neq j} C(\cT_{i,\infty})\bigg] - \sum_{j\geq 1} C(\cT_{i,\infty}) = \min_{j\geq 1}\bigg[W(\rho,j) - X_j\bigg]. \]
    This is exactly the same recursion as in (a) of the Theorem leading to the logistic distribution. Further if $j_{opt}$ is the child of $\rho$ that minimizes the right hand side then the matching of the root corresponds to the edge $(\rho, j_{opt})$. In \cite{aldous2001zeta}  Aldous found a magical scheme to make sense of  the above ``calculations'' and give a complete proof for the asymptotics of the limit constant. 
\end{enumeratea}

\subsection{Recurrence of random walks on planar graphs}
\begin{figure}[!htb]
    \centering
    \begin{minipage}{.5\textwidth}
        \centering
        \includegraphics[width=0.94\linewidth, height=0.35\textheight]{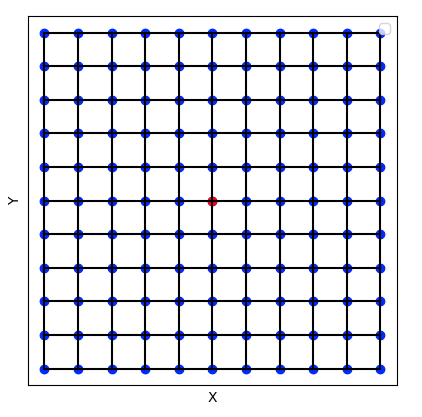}
        \caption{Standard integer lattice with the origin marked red. }
        \label{fig:lattice}
    \end{minipage}%
    \begin{minipage}{0.5\textwidth}
        \centering
        \includegraphics[width=0.94\linewidth, height=0.35\textheight]{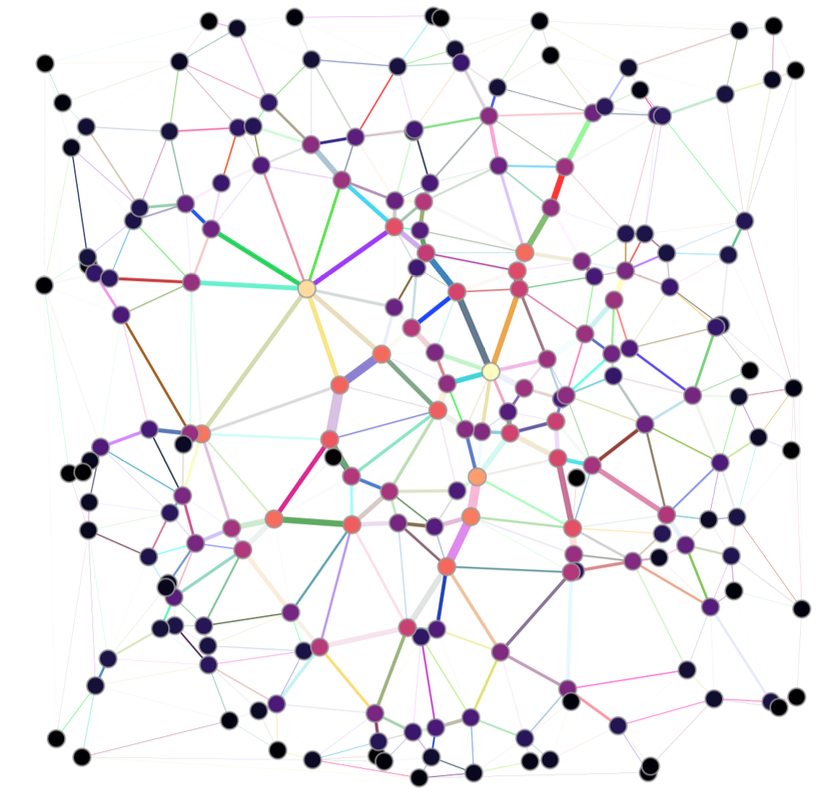}
        \caption{Random Triangulation of part of the plane resulting in a planar graph where each edge has degree $3$.}
        \label{fig:triang}
    \end{minipage}
\end{figure}

The previous section described the importance of ``spatial versions'' of local weak convergence, where edge lengths played an important role. We conclude this section with a remarkable result of Benjamini and Schramm in \cite{benjamini-schramm} that in fact was another genesis for the notion of local weak convergence in the early 2000s. To set the stage, recall that one of the classical results says that, if one considers the two dimensional integer lattice then the random walk is reccurrent, namely starting from any location, it returns to that location eventually with probability one. A natural question is: what if one considers other planar geometries, namely infinite graphs deterministic or random that can be embedded in space. Examples of two such geometries in (part of) the plane are given in Figures \ref{fig:lattice}-\ref{fig:triang}. We will now paraphrase the fundamental result in \cite{benjamini-schramm}, the proof idea is significantly different from the rest of this review so we refer the interested reader to this beautiful paper for details regarding the proof.

\begin{thm}[{\cite{benjamini-schramm}}]
    Let $\set{\cG_n:n\geq 1}$ be a sequence of (potentially random) finite planar connected graphs with maximal degrees uniformly bounded by a constant $M$. Suppose $\cG_n \LWC \cG_\infty^\star$ as $n\to\infty$. Then with probability one, the random walk on $\cG_\infty^\star$ is recurrent. 
\end{thm}
Once again what is amazing is that recurrence of the limit object which is in principle a completely global property is modulated by local weak convergence. 

\section{Local weak convergence for dynamic network models}
\label{sec:dynamic-res}

The previous Sections largely dealt with network models on \emph{static} networks, namely where the underlying topology of the network does not change. The goal here is to describe dynamic models, and for the sake of clarity we will focus on various families of growing random tree models, some of whose motivation was given in Section \ref{sec:mot} and fundamentals related to fringe convergence (Section \ref{sec:fringe}) were described in Section \ref{sec:growing-trees-fringe}. To keep this paper to manageable length, we will give a proof idea only in the first subsection and only paraphrase results in the remaining sections giving references to where the interested reader can find more details. 

\subsection{Local weak convergence and PageRank }
Recall PageRank from Definition \ref{def:PageRank}. Introduced by Google to rank web pages \cite{page1999pagerank}, PageRank is arguably one of the most effective centrality measures on complex networks and is a main factor behind the initial success of Google as a search engine. When compared to simpler centrality measures like degree, PageRank incorporates the effect of geometry of the network around a vertex beyond finite radius neighborhoods, and this should likely make it a more sensitive notion of centrality. This observation is elucidated by the following representation of PageRank, which we describe for tree networks for simplicity (an analogous representation holds for non-tree networks). For a tree network $\cG_n$ with $n$ vertices (viewed as a directed graph with edges pointing from offspring to parent), and $v\in \cG_n$, let $P_l(v,n)$ denote the number of \emph{directed} paths of length $l$ that end at $v$ in $\cG_n$. It is easy to check from Definition \ref{def:PageRank} that the PageRank scores have the explicit formulae for any vertex $v$, 
\begin{equation}
\label{eqn:non-root-page-rank}
\fR_{v,c}=\frac{(1-c)}{n}\left(1 + \sum_{l=1}^\infty c^l P_l(v,n)\right).
\end{equation}
For the sequel, it will be easier to formulate results in terms of the \emph{graph normalized} PageRank scores $\set{R_{v,c}:v\in \cG_n} = \set{n\fR_{v,c}(n): v\in \cG_n}$.

A natural question that one could ask is whether PageRank and degree truly differ at the \emph{extremal} (or large deviations) level. Namely, how do the identities of the most central vertices in the network quantified via degree and PageRank compare?  
Along this line of enquiry, a popular belief is the so-called \emph{power-law hypothesis}, which conjectures that for real world networks with a power-law degree distribution, the PageRank also has a power-law distribution with the same exponent as the degree. This has been shown to hold in several \emph{static} network models (no time evolution) like the directed configuration model \cite{chen2017generalized, olvera2019pagerank} and the inhomogeneous random digraph \cite{lee2020pagerank, olvera2019pagerank}. For dynamic random graphs, the picture was recently shown to be very different. For linear preferential attachment models, \cite{banerjee2021pagerank} showed that the PageRank distribution has heavier tails than degree. This phenomenon has been extended to a variety of other dynamic random graph models in \cite{banerjee2022co,antunes2023attribute}. Strikingly, for the random recursive tree described in Section \ref{sec:growing-trees-fringe}, the limiting degree distribution has an exponentially decaying tail while the limiting PageRank distribution has a power-law tail with exponent $1/c$.

To prove such results, local weak convergence serves as a crucial tool. At a high level, the idea is sketched as follows.

(a) \textbf{Enter local weak convergence: } Observe from \eqref{eqn:non-root-page-rank} that the PageRank of a vertex $v$ depends only on its \emph{in-component}, namely the subgraph spanned by vertices with a directed path to $v$. Moreover, although the PageRank of $v$ depends on the whole in-component of $v$, the contribution of a vertex $u$ to the PageRank of $v$ decays exponentially with its distance from $v$. From these observations, it is not too hard to believe that if the in-component of a uniformly chosen vertex in the network converges in a (directed) local weak sense to a rooted limiting random tree $\cG_{\infty}$, then the limiting PageRank distribution (limit of the empirical distribution of PageRanks of vertices in $\cG_n$ as $n \rightarrow \infty$) corresponds to the law of the root PageRank in $\cG_{\infty}$. This was rigorously verified in \cite{garavaglia2020local,banerjee2021pagerank}. 

(b) \textbf{PageRank on static graphs: } The distribution of the root PageRank in $\cG_{\infty}$ is qualitatively different for static and dynamic random graphs. For static graphs, the (directed) local weak limit corresponds to the progeny tree of a Galton-Watson Branching process. This is not a surprise, given our discussion in Section \ref{staticlwl}. This implies that the root PageRank $\mathcal{R}_{\emptyset}$ in the limiting tree satisfies the \emph{recursive distributional fixed point equation}
   \begin{equation}\label{rde}
   \mathcal{R}_{\emptyset} \, \equald \, c\sum_{i=1}^{\mathcal{D}}\mathcal{R}_i + (1-c),
   \end{equation}
   where $\mathcal{D}$ is the (in-)degree of the root and $\{\mathcal{R}_i\}$ are iid having the same distribution as $\mathcal{R}$. This equation can then be analyzed via renewal theory to establish that 
   $$
   0< \liminf_{x \rightarrow \infty}\frac{\prob[\mathcal{R}_{\emptyset}>x]}{ \prob[\mathcal{D}>x]} \le \limsup_{x \rightarrow \infty}\frac{\prob[\mathcal{R}_{\emptyset}>x]}{ \prob[\mathcal{D}>x]} < \infty.
   $$
   A detailed proof can be found in \cite{jelenkovic2010information}. This line of argument verifies the power law hypothesis for a variety of static random graphs.

 (c) \textbf{PageRank on dynamic graphs: }For dynamic graphs, as previously indicated, the local weak limits turn out to be more exotic.
 For a concrete example, consider the \emph{non-uniform random recursive tree} introduced in Section \ref{sec:growing-trees-fringe}. As previously discussed, the discrete network process can be embedded in a \emph{continuous time} branching processes $\{\BP(t): t \ge 0\}$ where each individual in the population reproduces (independently across individuals) at times given by a point process $\zeta(\cdot)$ on $\mathbb{R}_+$ with intensity measure $\mu(\cdot)$ and newly born individuals start reproducing according to an independent copy of the same point process. In the point process $\zeta(\cdot)$, the time gap between the $k$-th and $(k+1)$-th birth has distribution $\exp(f(k))$, where $f$ is the attachment function governing the network evolution. The directed local weak limit is then given by $\BP(T_{\lambda})$, where $T_{\lambda} \sim \exp(\lambda)$ is independent of $\BP(\cdot)$ and $\lambda$ is the \emph{Malthusian rate} of growth of the branching process, given analytically (under some regularity assumptions) as the unique root of the equation
  $$
  \int_0^{\infty} e^{-\lambda t}\mu(dt) = 1.
  $$
  See \cite{rudas2007random} for a detailed treatment.
  Observe that different vertices are born at different times and hence, at any time $t$, the distribution of the number of offspring of each existing child of the root differ across children. It is this \emph{temporal inhomogeneity} that is behind the fundamentally different behavior of PageRank in dynamic networks as one cannot exploit a recursive distributional equation analogous to \eqref{rde}. However, there is a different way to look at things in terms of \emph{percolation} on branching processes described below.

  \begin{defn}[Percolation on $\BP$]
	\label{def:perc-bp}
	Fix a damping factor $c\in (0,1)$. For any $t\geq 0$, write $\BP^c(t)$ for the connected cluster of the root (which is also a tree) when we retain each edge $e\in \BP(t)$ with probability $c$ and delete with probability $(1-c)$, independently across edges. Write $\set{\BP^c(t):t\geq 0}$ for the corresponding non-decreasing rooted tree valued Markov process where children born to vertices in the connected cluster of the root are retained with probability $c$ at their time of birth. 
\end{defn} 
By \eqref{eqn:non-root-page-rank}, it then follows that
\begin{equation}
\label{eqn:217}
	\mathcal{R}_{\emptyset} \, \equald \, (1-c)\E\left(|\BP^c(T_{\lambda})| \ \Big| \ \BP(T_{\lambda})\right). 
\end{equation}
It turns out that $\BP^c(\cdot)$ is a continuous time branching process in its own right with Malthusian rate $\lambda^c$ (say). Then, heuristically, recalling that $|\BP(t)| \sim e^{\lambda^c t}$, as $x \rightarrow \infty$,
$$
\Pr[\mathcal{R}_{\emptyset} > x] \sim \Pr[e^{\lambda^c T_{\lambda}} > x] \sim \Pr\left[T_{\lambda} > \frac{1}{\lambda^c}\log x\right] = x^{-\lambda/\lambda^c}.
$$
This idea has been rigorously laid out in a number of cases. We present a result from \cite{banerjee2021pagerank} summarizing this for the linear preferential attachment case.

\begin{theorem}\label{exponent}
Consider the linear preferential attachment tree process $\{\cG_n : n \ge 0\}$ with attachment function $f(k) = k
+1+\beta, \, k \ge 0,$ for some fixed $\beta>0$. The in-degree and PageRank of a uniformly chosen vertex $V_n$ in $\cG_n$ jointly converge in distribution:
\begin{equation*}\label{weakcon0}
(D_{V_n}, R_{V_n}) \xrightarrow{d} (\mathcal{D}, \mathcal{R}_\emptyset) \ \text{ as } n \rightarrow \infty.
\end{equation*}
Moreover, there exists $C>0$ such that as $k \rightarrow \infty$,
\begin{equation*}\label{degtail}
\Pr[\mathcal{D} \ge k] = Ck^{-2-\beta}\left(1 + O(k^{-1})\right),
\end{equation*}
and positive constants $C_1, C_2$ such that for any $r \geq 1$, 
\begin{equation*}\label{prtail}
C_1r^{-(2+\beta)/(1+(1+\beta)c)} \ \le \ \Pr[\mathcal{R}_{\emptyset} \ge r] \ \le \ C_2r^{-(2+\beta)/(1+(1+\beta)c)}.
\end{equation*}
\end{theorem}
More examples are discussed in the next few Subsections. 

We remark here that, although the above discussion was presented for tree networks, the ideas extend naturally to \emph{locally tree-like} networks. For such networks, the local weak limits are trees and the arguments sketched in (b) and (c) above apply to the limiting objects. Local weak limits for non-tree linear preferential attachment models are described in \cite{berger2014asymptotic,banerjee2021pagerank,garavaglia2022universality}, and limits for collapsed branching processes are obtained in \cite{banerjee2023local}.

\subsection{Attributed network models}

Next recall the nodal attribute model defined in Section \ref{sec:nodal-att-mot}. The following paraphrases some of the main results in \cite{antunes2023attribute}. 

\begin{thm}[{\cite{antunes2023attribute}}]
    Consider the linear ($\gamma\equiv 1$) setting of the model with finite attribute space and assume the propensity kernel $\kappa > 0$ (i.e. every entry of this matrix is strictly positive). Then:
    \begin{enumeratea}
        \item While the model cannot be embedded directly as a continuous time, multitype, branching process, its evolution can be analyzed using stochastic approximation techniques to show that the sequence of trees converges to a limiting infinite {\tt sin}-tree with nodal types. 
        \item This result gives information on joint  distribution asymptotics of types and degrees, showing that tail exponents of the limiting degree distribution of different types can depend on the type. 
        \item The asymptotic limit also gives information about the PageRank distribution, in particular showing that the extremal behavior of the PageRank scores do not depend on types. 
        \item These limit results and constructions related to these results also give information on the behavior of various network sampling algorithms. In particular, it is shown in some settings with rare minority vertices that sampling from the graph with probability proportional to PageRank scores, which can be accomplished using local exploration schemes, has a quantifiably higher chance of sampling minority vertices in comparison to uniform and degree-based sampling schemes. 
    \end{enumeratea}
\end{thm}

\subsection{Co-evolving networks}
 One major frontier, especially for developing rigorous understanding of proposed models, are the so-called co-evolutionary (or adaptive) networks, where specific dynamics (e.g. random walk explorations) on the network influence the structure of the network, which in turn influences the dynamics; thus both modalities (dynamics on the network and the network itself) co-evolve \cites{gross2008adaptive,aoki2016temporal,sayama2013modeling,sayama2015social}. Motivated by the growth of social networks, there has been significant interest in trying to understand the influence of processes such as search engines or influence ranking mechanisms in the growth of networks. A number of papers \cites{pandurangan2002using,blum2006random,chebolu2008pagerank} have explored the dynamic evolution of networks through new nodes first exploring neighborhoods of randomly selected vertices before deciding on whom to connect.  We now describe a specific class of such models.

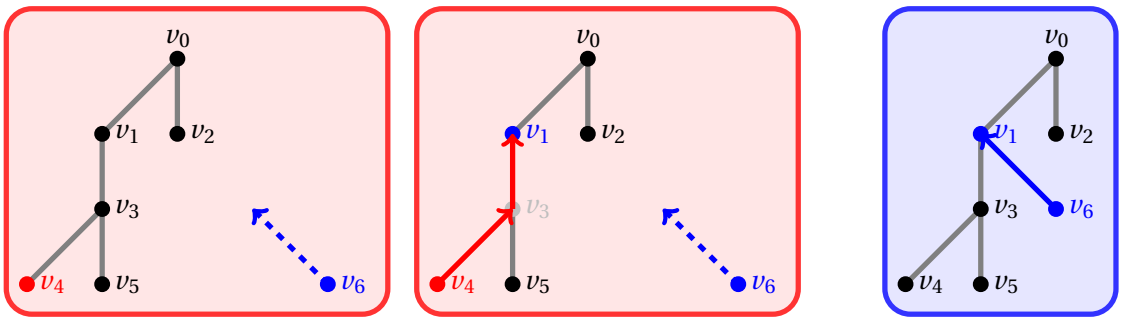
\begin{figure}[h]
     \centering     
\begin{subfigure}[b]{0.3\textwidth}
         \centering
         \begin{center}
\begin{tikzpicture}[background rectangle/.style=
		     {draw=red!80,fill=red!10,rounded corners=2ex},
		   show background rectangle]
  \draw[gray](0,0)--(0,-1);
  \draw[gray](0,0)--(-1,-1);
  \draw[gray](-1,-1)--(-1,-2);
  \draw[gray](-1,-2)--(-1,-3);
  \draw[gray](-1,-2)--(-2,-3);
    \filldraw(0,0) circle (2pt) node[above] {$v_0$};
   \filldraw(0,-1) circle (2pt) node[right] {$v_2$};
    \filldraw(-1,-1) circle (2pt) node[right] {$v_1$};
    \filldraw(-1,-2) circle (2pt) node[right] {$v_3$};
     \filldraw(-1,-3) circle (2pt) node[right] {$v_5$};
       \filldraw[red](-2,-3) circle (2pt) node[right] {$v_4$};
       \draw[dashed, blue,->](2,-3)--(1,-2);
       \filldraw[blue](2,-3) circle (2pt) node[right] {$v_6$};
\end{tikzpicture}
\end{center}
\end{subfigure}
\hfill
\begin{subfigure}[b]{0.3\textwidth}
         \centering
         \begin{center}
\begin{tikzpicture}[background rectangle/.style=
		     {draw=red!80,fill=red!10,rounded corners=2ex},
		   show background rectangle]
  \draw[gray](0,0)--(0,-1);
  \draw[gray](0,0)--(-1,-1);
  \draw[gray](-1,-2)--(-1,-3);
 \filldraw(0,0) circle (2pt) node[above] {$v_0$};
   \filldraw(0,-1) circle (2pt) node[right] {$v_2$};
    \filldraw[blue](-1,-1) circle (2pt) node[right] {$v_1$};
    \filldraw[lightgray](-1,-2) circle (2pt) node[right] {$v_3$};
     \filldraw(-1,-3) circle (2pt) node[right] {$v_5$};
       \filldraw[red](-2,-3) circle (2pt) node[right] {$v_4$};
         \draw[->, red](-1,-2)--(-1,-1);
          \draw[->, red](-2,-3)--(-1,-2);   
          
              \draw[dashed, blue,->](2,-3)--(1,-2);
       \filldraw[blue](2,-3) circle (2pt) node[right] {$v_6$};
\end{tikzpicture}
\end{center}
     \end{subfigure}
\hfill     
\begin{subfigure}[b]{0.3\textwidth}
         \centering
     \begin{center}
\begin{tikzpicture}[background rectangle/.style=
		     {draw=blue!80,fill=blue!10,rounded corners=2ex},
		   show background rectangle]
  \draw[gray](0,0)--(0,-1);
  \draw[gray](0,0)--(-1,-1);
  \draw[gray](-1,-1)--(-1,-2);
  \draw[gray](-1,-2)--(-1,-3);
  \draw[gray](-1,-2)--(-2,-3);
    \filldraw(0,0) circle (2pt) node[above] {$v_0$};
   \filldraw(0,-1) circle (2pt) node[right] {$v_2$};
    \filldraw[blue](-1,-1) circle (2pt) node[right] {$v_1$};
    \filldraw(-1,-2) circle (2pt) node[right] {$v_3$};
     \filldraw(-1,-3) circle (2pt) node[right] {$v_5$};
       \filldraw(-2,-3) circle (2pt) node[right] {$v_4$};
          \draw[->, blue](0,-2)--(-1,-1);
        \filldraw[blue](0,-2) circle (2pt) node[right] {$v_6$};
\end{tikzpicture}
\end{center}
     \end{subfigure}
        \caption{$v_6$ is a new incoming vertex, and selects $v_4$ to start exploring the network,  with sampled number of exploration steps $Z_6=2$.}
        \label{fig:three graphs}
\end{figure}

Fix a probability mass function $\vp:= \set{p_k:k\geq 0}$ on $\bZ_+$. For the rest of the paper, let $\vZ = \set{Z_1, Z_2, \ldots}$ be an \emph{i.i.d} sequence with distribution $\vp$.  We now describe the recursive construction of a sequence of random trees $\set{\cT_n:n\geq 1}$, always rooted at vertex $\set{v_0}$, with edges pointed from descendants to their parents.  Start with two vertices $\set{v_0, v_1}$, with $\cT_1$ a rooted tree at $\set{v_0}$, an oriented edge from $v_1$ to $v_0$.  Assume for some $n\geq 1$, we have constructed $\cT_n$. Then to construct $\cT_{n+1}$:

\begin{enumeratea}
    \item New vertex $\set{v_{n+1}}$ enters the system at time $n+1$. 
    \item This new vertex selects a vertex $V_n$, \emph{uniformly at random}, amongst the existing vertices $\cV(\cT_n) = \set{v_0, \ldots, v_n} $. 
    \item  Let $\cP(v_0, V_n)$ denote the path from the root to this vertex. This new vertex traverses up this path for a random length $Z_{n+1}$ and attaches to the terminal vertex. If the {graph distance to the root,} $\dist(v_0, V_n) \leq Z_{n+1}$ then this new vertex attaches to the root $v_0$. 
\end{enumeratea}

It turns out, the setting where the pmf $\vp$ satisfies $p_0+p_1 =1$ is identical to the linear preferential attachment model. So to state our main results, we will assume that $p_0+p_1 <1$ and further $\E(Z) < \infty$. To get some intuition consider the simulation figures below (network size $n=30,000$) and $\vZ$ having a Geometric distribution with parameter $p=.35, .65$ respectively. 

\begin{figure}[!htb]
    \centering
    \begin{minipage}{.5\textwidth}
        \centering
        \includegraphics[width=0.94\linewidth, height=0.35\textheight]{coevoling_simulation_35.pdf}
        \caption{$\vZ\sim$ Geometric distribution with $p=.35$ }
        \label{fig:geom1}
    \end{minipage}%
    \begin{minipage}{0.5\textwidth}
        \centering
        \includegraphics[width=0.94\linewidth, height=0.35\textheight]{coevoling_simulation_65.pdf}
        \caption{$\vZ\sim$ Geometric distribution with $p=.65$}
        \label{fig:geom2}
    \end{minipage}
\end{figure}

When $p=.35$, this implies $\E(Z) > 1$ namely the exploration walk tends to explore further up the path to the root, whereas when $p=.65$ then $\E(Z)\leq 1$, the walk is more local in nature. It turns out $\E(Z)=1$ is the ``phase transition'' point; we paraphrase a few results from \cite{banerjee2022co}. 

\begin{thm}[\cite{banerjee2022co}]
    For the above network model:
    \begin{enumeratea}
        \item When $\E(Z) \leq 1$ then the sequence of random trees converges to a limiting infinite random sin-tree $\cT_n \probcrf \cT_\infty$ from which one can derive information about various quantities including asymptotics for the degree distribution and the corresponding degree exponents. This has deep connections with large deviations and quasi-stationary distributions of random walks on $\mathbb{Z}$. 
        \item When $\E(Z) >1$ then $\set{\cT_n:n\geq 1}$ converges in the fringe sense (Def \ref{def:local-weak} \eqref{it:fringe-a}) but not in the extended fringe sense. Further there is condensation at the root in the sense that there is a model dependent strictly positive constant $\gamma(\vp) >0$ such that the root degree $\deg(\rho,n)$ in $\cT_n$ satisfies,   $\deg(\rho,n)/n \probc \gamma(\vp)$ as $n \rightarrow \infty$.
        \item We noted before that the power-law hypothesis (comparing PageRank and degree tail behavior) holds for several static graphs and fails to hold for several dynamic graphs. The model under consideration `interpolates' these two regimes in the following sense. When $\E[Z]\leq 1$, then (under some regularity assumptions), there exist $R\ge 1$ and $s(\vp) \in (0,1)$ such that the limiting degree and PageRank distributions satisfy:
        $$
        \lim_{k\to\infty} \frac{ \log(\pr(\cD\geq k))}{\log k} = -R,
        $$
$$\lim_{r\to\infty} \frac{ \log(\pr(\cR_{\emptyset}\geq r))}{\log r}=
\begin{cases}
-R &\text{ for $c\in(0,s(\vp)]$ with $c<1$},\\
-\frac{1}{cf(1/c)} & \text{ for $c\in (s(\vp),1)$}.
\end{cases}
$$
    \end{enumeratea} 
    Hence, by varying the damping factor, one transitions from the regime where the power-law hypothesis holds to one where PageRank has strictly heavier tails than the degree distribution. Moreover, one can verify that $\frac{1}{s(\vp)f(1/s(\vp))} = R$, and thus this transition is continuous.
\end{thm}


\section{Conclusion and further reading}
\label{sec:conc}
As mentioned in the introduction, this paper is not meant to be a survey and rather was meant to be a gentle introduction to junior researchers in the use of local weak convergence in various problems of probabilistic combinatorics; still we apologize to the many wonderful researchers who have contributed to this methodology and its use at the frontiers of the field that we were not able to cite. Our goal now is to give starting points for such readers to take next steps in their exploration of this vast field and get closer to the edge of research and the vast unexplored terrain of the unknown. 

For such researchers, perhaps the most comprehensive survey can be found in \cite[Chapter 2]{van2023random} which outlines all the myriad fundamentals and their extensions, especially in the context of unweighted graphs,  in great detail. In the geometric setting, for the fundamentals and their use in probabilistic combinatorial optimization, \cite{aldous-steele-obj} is a phenomenal resource. We also urge readers to peruse \cite{aldous-fringe} which laid the genesis for many of these concepts via considering the setting of trees and which also explains and connects these topics to the more classical findings related to stable age distribution theory of continuous time branching processes by Jagers and Nerman. 

For statistical physics inspired models such as the Ising model on sparse random graphs, the two summer school lecture notes \cite{10.1214/09-BJPS027,van2017stochastic} are wonderful resources for rigorous theory, while \cite{mezard2009information} gives an overview of the intuition, especially from statistical physics for what ``should happen''. One major tool that seemed to show up repeatedly in Section \ref{sec:power} and Section \ref{sec:dynamic-res} were recursive distributional equations and \cite{10.1214/105051605000000142} gives a wide-ranging overview of rigorous theory of such constructs. 

Finally for random matrix theory and its connections to local weak convergence, the various lecture notes of Charles Bordenave e.g. \cite{bordenave2016spectrum} are a great starting point, while for planar graphs, the St. Flour lecture notes of Asaf Nachmias \cite{nachmias2020planar} provide a great diving board into the ramifications of one of the classics in this entire field, namely the work of Benjamini and Schramm \cite{benjamini-schramm}.

\section*{Acknowledgement}
Banerjee was supported in part by the NSF CAREER award DMS-2141621. Bhamidi was supported in part by NSF DMS-2113662. Banerjee, Bhamidi and Young were partially funded by NSF RTG grant DMS-2134107 .  We thank David Aldous, Remco van der Hofstad, Mariana Olvera-Cravioto and Allan Sly for many insights over the years.  

\bibliographystyle{plain}

\end{document}